# A UNIFIED APPROACH TO MODEL SELECTION AND SPARSE RECOVERY USING REGULARIZED LEAST SQUARES[1]

By Jinchi Lv and Yingying Fan

*University of Southern California*

Model selection and sparse recovery are two important problems for which many regularization methods have been proposed. We study the properties of regularization methods in both problems under the unified framework of regularized least squares with concave penalties. For model selection, we establish conditions under which a regularized least squares estimator enjoys a nonasymptotic property, called the weak oracle property, where the dimensionality can grow exponentially with sample size. For sparse recovery, we present a sufficient condition that ensures the recoverability of the sparsest solution. In particular, we approach both problems by considering a family of penalties that give a smooth homotopy between $L_0$ and $L_1$ penalties. We also propose the sequentially and iteratively reweighted squares (SIRS) algorithm for sparse recovery. Numerical studies support our theoretical results and demonstrate the advantage of our new methods for model selection and sparse recovery.

**1. Introduction.** Model selection and sparse recovery are two important areas that have attracted much attention of the researchers. They are different but related, and share some common ideas especially when dealing with large scale problems. Examples include the lasso [Tibshirani (1996)], SCAD [Fan (1997) and Fan and Li (2001)], Dantzig selector [Candes and Tao (2007)] and MCP [Zhang (2007)] in model selection, and the basis pursuit [Chen, Donoho and Saunders (1999)] and many other $L_1$ methods [Candes and Tao (2005, 2006) and Candès, Wakin and Boyd (2008)] in sparse recovery. The analysis of vast data sets with the number of variables $p$ comparable to or much larger than the number of observations $n$ frequently arises nowadays

Received February 2008; revised January 2009.

[1]Supported in part by NSF Grants DMS-08-06030 and DMS-09-06784 and 2008 Zumberge Individual Award from USC's James H. Zumberge Faculty Research and Innovation Fund.

in both areas and poses many challenges that are not present in smaller scale studies. Sparsity plays an important role in these large scale problems. It is often believed that only a small fraction of the data is informative, whereas most of it is noise.

Consider the linear regression model

$$\mathbf{y} = \mathbf{X}\boldsymbol{\beta} + \boldsymbol{\varepsilon}, \tag{1}$$

where $\mathbf{y}$ is an $n$-dimensional vector of responses, $\mathbf{X} = (\mathbf{x}_1, \ldots, \mathbf{x}_p)$ is an $n \times p$ design matrix, $\boldsymbol{\beta} = (\beta_1, \ldots, \beta_p)^T$ is an unknown $p$-dimensional vector of regression coefficients and $\boldsymbol{\varepsilon}$ is an $n$-dimensional vector of noises. In the sparse modeling, we assume that a fraction of the true regression coefficients vector $\boldsymbol{\beta}_0 = (\beta_{0,1}, \ldots, \beta_{0,p})^T$ are exactly zero. In this paper we allow $\boldsymbol{\beta}_0$ to depend on $n$. We denote by $\mathfrak{M}_0 = \mathrm{supp}(\boldsymbol{\beta}_0)$ the support of $\boldsymbol{\beta}_0$, which is called the true underlying sparse model hereafter. Model selection aims to locate those predictors $\mathbf{x}_j$ with nonzero $\beta_{0,j}$, which are called true variables hereafter, and to give consistent estimate of $\boldsymbol{\beta}_0$ on its support. Throughout the paper we consider deterministic design matrix $\mathbf{X}$ and assume the identifiability of $\boldsymbol{\beta}_0$, in the sense that the equation $\mathbf{X}\boldsymbol{\beta}_0 = \mathbf{X}\boldsymbol{\beta}$, $\boldsymbol{\beta} \in \mathbf{R}^p$ entails either $\boldsymbol{\beta} = \boldsymbol{\beta}_0$ or $\|\boldsymbol{\beta}\|_0 > \|\boldsymbol{\beta}_0\|_0$. We denote by $\|\cdot\|_q$ the $L_q$ norm on the Euclidean spaces $q \in [0, \infty]$. Many methods have been proposed in the literature to construct estimators that mimic the oracle estimators under different losses, where the oracle knew the true model $\mathfrak{M}_0$ ahead of time. The main difficulty of recovering $\mathfrak{M}_0$ lies in the collinearity among the predictors, which increases as the dimensionality grows. See, for example, Fan and Li (2006) for a comprehensive overview of challenges of high dimensionality in statistics.

Much insight into model selection can be obtained if we understand a closely related problem of sparse recovery, which aims at finding the minimum $L_0$ (sparsest possible) solution to the linear equation

$$\mathbf{y} = \mathbf{X}\boldsymbol{\beta}, \tag{2}$$

where $\boldsymbol{\beta} = (\beta_1, \ldots, \beta_p)^T$, $\mathbf{y} = \mathbf{X}\boldsymbol{\beta}_0$ and $\mathbf{X}$ and $\boldsymbol{\beta}_0$ are the same as those in the linear model (1). The identifiability of $\boldsymbol{\beta}_0$ assumed above ensures that our target solution here is unique and exactly $\boldsymbol{\beta}_0$. When the $p \times p$ matrix $\mathbf{X}^T\mathbf{X}$ is singular or close to singular, finding $\boldsymbol{\beta}_0$ is not an easy task. It is known that directly solving the $L_0$-regularization problem is combinatorial and, thus, is impractical in high dimensions. To attenuate this difficulty, many regularization methods such as the $L_1$ method of basis pursuit have been proposed to recover $\boldsymbol{\beta}_0$, where continuous penalty functions are used in place of the $L_0$ penalty. This raises a natural question: under what conditions does a regularization method give the same solution as that of the $L_0$ regularization? Many authors have contributed to identifying conditions that ensure the $L_1/L_0$ equivalence. In this paper we generalize a sufficient condition identified for the $L_1$ penalty to concave penalties, which ensures such an equivalence.

In view of (1) and (2), model selection and sparse recovery can be regarded as two interrelated problems. Due to the presence of noise, recovering the true model $\mathfrak{M}_0$ in (1) is intrinsically more challenging than recovering the sparsest possible solution $\boldsymbol{\beta}_0$ in (2). Various regularization methods in model selection such as those mentioned before have been studied by many researchers. See, for example, Bickel and Li (2006) for a comprehensive review of regularization methods in statistics.

In a seminal paper, Fan and Li (2001) lay down the theoretical foundation of nonconvex penalized least squares and nonconcave penalized likelihood for variable selection, and introduce the concept of model selection oracle property. An estimator $\widehat{\boldsymbol{\beta}}$ is said to have the oracle property [Fan and Li (2001)] if: (1) it enjoys the sparsity in the sense of $\widehat{\boldsymbol{\beta}}_{\mathfrak{M}_0^c} = \mathbf{0}$ with probability tending to 1 as $n \to \infty$, and (2) it attains an information bound mimicking that of the oracle estimator [see also Donoho and Johnstone (1994)], where $\widehat{\boldsymbol{\beta}}_{\mathfrak{M}_0^c}$ is a subvector of $\widehat{\boldsymbol{\beta}}$ formed by components with indices in $\mathfrak{M}_0^c$, the complement of the true model $\mathfrak{M}_0$. Fan and Li (2001) study the oracle properties of nonconcave penalized likelihood estimators in the finite-dimensional setting. Their results were extended later by Fan and Peng (2004) to the setting of $p = o(n^{1/5})$ or $o(n^{1/3})$. In this paper we generalize the results of Fan and Li (2001) and Fan and Peng (2004), in the setting of regularized least squares, to a more general triple $(s, n, p)$ for concave penalties, where $s = \|\boldsymbol{\beta}_0\|_0$. In particular, we show that under some regularity conditions, the regularized least squares estimator enjoys a nonasymptotic weak oracle property, where the dimensionality $p$ can be of exponential order in sample size $n$. This constitutes one of the main contributions of the paper.

In this paper, we consider both problems of model selection and sparse recovery in the unified framework of regularized least squares with concave penalties. Specifically, for sparse recovery we construct the solutions of regularization problems under the constraint in (2) by analyzing the solutions of related regularized least squares problems and then letting the regularization parameter $\lambda \to 0+$. In particular, we consider a family of penalty functions that give a smooth homotopy between $L_0$ and $L_1$ penalties for both problems. The unified approach using the $L_1$ penalty has been considered by Chen, Donoho and Saunders (1999), Fuchs (2004), Donoho, Elad and Temlyakov (2006) and Tropp (2006), among others.

The rest of the paper is organized as follows. In Section 2, we discuss the choice of penalty functions. We study the properties of regularization methods in model selection and sparse recovery for concave penalties in Sections 3 and 4. Section 5 discusses algorithms for solving regularization problems. In Section 6 we present four numerical examples using both simulated and real data sets. Proofs are presented in Section 7. We provide some discussion of our results and their implications in Section 8.

**2. Regularization methods with concave penalties.** In this paper, we study regularization methods in model selection and sparse recovery for concave penalties. For sparse recovery in (2), we consider the $\rho$-regularization problem

$$\text{(3)} \qquad \min \sum_{j=1}^{p} \rho(|\beta_j|) \quad \text{subject to} \quad \mathbf{y} = \mathbf{X}\boldsymbol{\beta},$$

where $\rho(\cdot)$ is a penalty function and $\boldsymbol{\beta} = (\beta_1, \ldots, \beta_p)^T$. For model selection in (1), we consider the regularized least squares problem

$$\text{(4)} \qquad \min_{\boldsymbol{\beta} \in \mathbf{R}^p} \left\{ 2^{-1} \|\mathbf{y} - \mathbf{X}\boldsymbol{\beta}\|_2^2 + \Lambda_n \sum_{j=1}^{p} p_{\lambda_n}(|\beta_j|) \right\},$$

where $\Lambda_n \in (0, \infty)$ is a scale parameter, $p_{\lambda_n}(\cdot)$ is a penalty function, $\lambda_n \in [0, \infty)$ is a regularization parameter indexed by sample size $n$ and $\boldsymbol{\beta} = (\beta_1, \ldots, \beta_p)^T$. We will drop the subscript $n$ when it causes no confusion. For any penalty function $p_\lambda$, let $\rho(t; \lambda) = \lambda^{-1} p_\lambda(t)$ for $t \in [0, \infty)$ and $\lambda \in (0, \infty)$. For simplicity, we will slightly abuse the notation and write $\rho(t; \lambda)$ as $\rho(t)$ when there is no confusion.

2.1. *Penalty functions.* By the nature of sparse recovery, the $L_0$ penalty $\rho(t) = I(t \neq 0)$ is the target penalty in (3), whereas other penalties may also be capable of recovering $\boldsymbol{\beta}_0$. As mentioned before, the $L_0$ penalty is not appealing from the computational point of view due to its discontinuity. It is known that the $L_2$ penalty $\rho(t) = t^2$ in (3) or (4) is analytically tractable, but generally produces nonsparse solutions. Such concerns have motivated the use of penalties that are computationally tractable approximations or relaxations to the $L_0$ penalty. Among all proposals the $L_1$ penalty $\rho(t) = t$, $t \in [0, \infty)$, has attracted much attention of the researchers in both sparse recovery and model selection. It has been recognized that the $L_1$ penalty is not an oracle that always points us to the true underlying sparse model.

Hereafter, we consider penalty functions $\rho(\cdot)$ that satisfy the following condition.

CONDITION 1. $\rho(t)$ is increasing and concave in $t \in [0, \infty)$, and has a continuous derivative $\rho'(t)$ with $\rho'(0+) \in (0, \infty)$. If $\rho(t)$ is dependent on $\lambda$, $\rho'(t; \lambda)$ is increasing in $\lambda \in (0, \infty)$ and $\rho'(0+)$ is independent of $\lambda$.

Fan and Li (2001) advocate penalty functions that give estimators with three desired properties—unbiasedness, sparsity and continuity—and provide insights into them [see also Antoniadis and Fan (2001)]. We discuss the connection of Condition 1 with these properties. Consider problem (4) with $n \times 1$ orthonormal design matrix and $\Lambda_n = 1$,

$$(5) \qquad \min_{\theta \in \mathbf{R}} \{2^{-1}(z - \theta)^2 + p_\lambda(|\theta|)\},$$

where $z = \mathbf{X}^T \mathbf{y}$ and $p_\lambda(t) = \lambda \rho(t)$, $t \in [0, \infty)$. We denote by $\widehat{\theta}(z)$ the minimizer of problem (5). Fan and Li (2001) demonstrate that for the resulting estimator $\widehat{\theta}(z)$: (1) unbiasedness requires that the derivative $p'_\lambda(t)$ is close to zero when $t \in [0, \infty)$ is large, (2) sparsity requires $p'_\lambda(0+) > 0$ and (3) continuity with respect to data $z$ requires that the function $t + p'_\lambda(t)$, $t \in [0, \infty)$ attains its minimum at $t = 0$. Note that the concavity of $\rho$ in Condition 1 entails that $\rho'(t)$ is decreasing in $t \in [0, \infty)$. Thus penalties satisfying Condition 1 and $\lim_{t \to \infty} \rho'(t) = 0$ enjoy the unbiasedness and sparsity. However, the continuity does not generally hold for all penalties in this class. The SCAD penalty [Fan and Li (2001)] $p_\lambda(t)$, $t \in [0, \infty)$, is given by

$$(6) \qquad p'_\lambda(t) = \lambda \left\{ I(t \leq \lambda) + \frac{(a\lambda - t)_+}{(a - 1)\lambda} I(t > \lambda) \right\} \qquad \text{for some } a > 2,$$

where often $a = 3.7$ is used, and MCP [Zhang (2007)] $p_\lambda(t)$, $t \in [0, \infty)$, is given by $p'_\lambda(t) = (a\lambda - t)_+/a$. Both SCAD and MCP with $a \geq 1$ satisfy Condition 1 and the above three properties simultaneously. Although the $L_1$ penalty satisfies Condition 1 as well as sparsity and continuity, it does not enjoy the unbiasedness, since its derivative is identically one regardless of $t \in [0, \infty)$.

For a penalty function $\rho$, we define its maximum concavity as

$$(7) \qquad \kappa(\rho) = \sup_{t_1, t_2 \in (0, \infty), t_1 < t_2} -\frac{\rho'(t_2) - \rho'(t_1)}{t_2 - t_1}$$

and we define the local concavity of the penalty $\rho$ at $\mathbf{b} = (b_1, \ldots, b_q)^T \in \mathbf{R}^q$ with $\|\mathbf{b}\|_0 = q$ as

$$(8) \qquad \kappa(\rho; \mathbf{b}) = \lim_{\epsilon \to 0+} \max_{1 \leq j \leq q} \sup_{t_1, t_2 \in (|b_j| - \epsilon, |b_j| + \epsilon), t_1 < t_2} -\frac{\rho'(t_2) - \rho'(t_1)}{t_2 - t_1}.$$

By the concavity of $\rho$ in Condition 1, we have $0 \leq \kappa(\rho; \mathbf{b}) \leq \kappa(\rho)$. It is easy to show by the mean-value theorem that $\kappa(\rho)$ defined in (7) equals $\sup_{t \in (0, \infty)} -\rho''(t)$ and $\kappa(\rho; \mathbf{b})$ defined in (8) equals $\max_{1 \leq j \leq q} -\rho''(|b_j|)$, provided that $\rho$ has a continuous second derivative $\rho''(t)$.

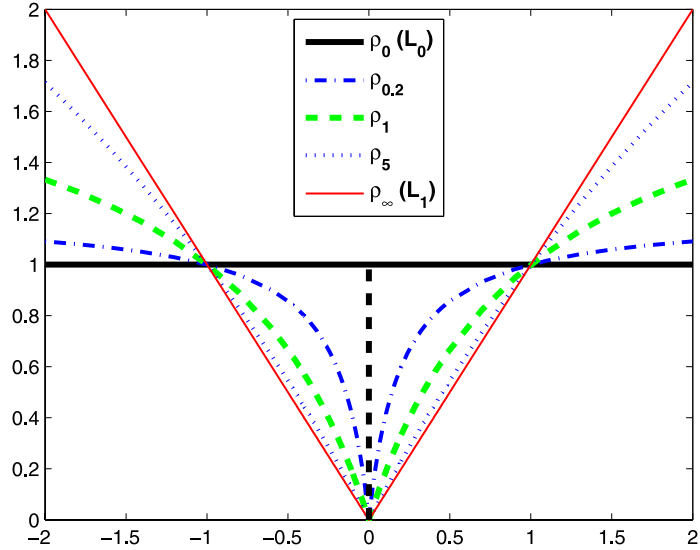

FIG. 1. *Plot of penalty functions $\rho_0$ ($L_0$, thick solid), $\rho_{0.2}$ (dash-dot), $\rho_1$ (dashed), $\rho_5$ (dotted), and $\rho_\infty$ ($L_1$, thin solid).*

2.2. *A family of penalties.* Nikolova (2000) studies the transformed $L_1$ penalty function $\rho(t) = bt/(1+bt)$, $t \in [0, \infty)$ and $b > 0$. A slight modification of it gives a family of penalties $\{\rho_a : a \in [0, \infty]\}$ given by, for $a \in (0, \infty)$,

$$(9) \qquad \rho_a(t) = \frac{(a+1)t}{a+t} = \left(\frac{t}{a+t}\right)I(t \neq 0) + \left(\frac{a}{a+t}\right)t, \qquad t \in [0, \infty),$$

and

$$(10) \qquad \rho_0(t) = \lim_{a \to 0+} \rho_a(t) = I(t \neq 0) \quad \text{and} \quad \rho_\infty(t) = \lim_{a \to \infty} \rho_a(t) = t$$

for $t \in [0, \infty)$. Figure 1 depicts $\rho_a$ penalties for a few $a$'s. We see from (9) that this modified family has the interpretation of a smooth homotopy between $L_0$ and $L_1$ penalties. So we refer to them as the smooth integration of counting and absolute deviation (SICA) penalties. It is easy to show that $\rho_a$ penalty with $a \in (0, \infty]$ satisfies Condition 1, and for each $a \in (0, \infty)$, $\lim_{t \to \infty} \rho'_a(t) = 0$. Thus $\rho_a$ with $a \in (0, \infty)$ gives estimators satisfying the unbiasedness and sparsity. As mentioned before, the continuity requires that the function $t + \lambda \rho'_a(t)$, $t \in [0, \infty)$, attains its minimum at $t = 0$; that is, $a \in [a_0, \infty]$, where $a_0 = \lambda + \sqrt{\lambda^2 + 2\lambda}$. Therefore, $\rho_a$ penalties $\rho_a$ with $a \in [a_0, \infty)$ satisfy Condition 1 and the above three properties simultaneously, and share the same spirit as SCAD and MCP. In addition, $\rho_a$ is infinitely differentiable on $[0, \infty)$ for each $a \in (0, \infty]$. The idea of linearly combining $L_0$ and $L_1$ penalties was investigated by Liu and Wu (2007).

For each $a \in (0, \infty)$, the $\rho_a$ penalty is closely related to the log penalty

$$\rho_{1,a}(t) = (a+1)\log(1+a^{-1}t), \qquad t \in [0, \infty).$$

In fact, the $L_1$ penalty is the first-order approximation to both $a(a+1)^{-1}\rho_a(t)$ and $a(a+1)^{-1}\rho_{1,a}(t)$, and always dominates them. Also, we have $\rho_a(t) = t\rho'_{1,a}(t)$. Clearly,

(11)
$$\rho'_a(t) = \frac{a(a+1)}{(a+t)^2}, \qquad t \in (0, \infty) \text{ for } a \in (0, \infty),$$

$$\rho'_a(0+) = 1 + a^{-1} \quad \text{for } a \in (0, \infty) \quad \text{and} \quad \rho'_\infty(t) = 1.$$

It is easy to see that the maximum concavity of $\rho_a$ penalty is

(12) $$\kappa(\rho_a) = \sup_{t \in (0,\infty)} -\rho''_a(t) = \sup_{t \in (0,\infty)} \frac{2a(a+1)}{(a+t)^3} = 2(a^{-1} + a^{-2}),$$

which is the maximum curvature of the curve $\rho_a$. Clearly, $\kappa(\rho_a)$ is decreasing in $a$, $\lim_{a \to 0+} \kappa(\rho_a) = \infty$, and $\lim_{a \to \infty} \kappa(\rho_a) = 0 = \kappa(\rho_\infty)$. Therefore, parameter $a$ controls the maximum concavity of $\rho_a$ and regulates where it stands between $L_0$ and $L_1$ penalties.

**3. Sparse recovery.** In this section we consider the $\rho$-regularization problem (3) for sparse recovery in (2). It is known that when the $n \times p$ matrix $\mathbf{X}$ is of full column rank $p$, (2) has a unique solution $\boldsymbol{\beta} = (\mathbf{X}^T\mathbf{X})^{-1}\mathbf{X}^T\mathbf{y}$. Otherwise, it has an infinite number of solutions, all of which form a $q$-dimensional linear subspace

(13) $$\mathcal{A} = \{\boldsymbol{\beta} \in \mathbf{R}^p : \mathbf{y} = \mathbf{X}\boldsymbol{\beta}\}$$

of $\mathbf{R}^p$ with $q = p - \text{rank}(\mathbf{X})$. Of interest is the nontrivial case of $q > 0$.

3.1. *Identifiability of $\beta_0$.* As mentioned in the Introduction, the minimum $L_0$ solution to (2) is

(14) $$\boldsymbol{\beta}_0 = \underset{\boldsymbol{\beta} \in \mathcal{A}}{\arg\min} \|\boldsymbol{\beta}\|_0.$$

Donoho and Elad (2003) introduce the concept of spark and show that the uniqueness of $\boldsymbol{\beta}_0$ can be characterized by the spark($\mathbf{X}$) of $\mathbf{X}$, where $\tau = \text{spark}(\mathbf{X})$ is defined as the smallest possible number such that there exists a subgroup of $\tau$ columns from the $n \times p$ matrix $\mathbf{X}$ that are linearly dependent. The spark of a matrix can be very different from its rank. For instance, the $n \times (n+1)$ matrix $[I_n \mathbf{e}_1]$ is of full rank $n$ and yet has spark equal to 2, where $\mathbf{e}_1 = (1, 0, \ldots, 0)^T$. In particular, they proved that any $\boldsymbol{\beta} \in \mathcal{A}$ with $\|\boldsymbol{\beta}\|_0 < \text{spark}(\mathbf{X})/2$ meets $\boldsymbol{\beta}_0$. Thus $\boldsymbol{\beta}_0$ is unique as long as $\|\boldsymbol{\beta}_0\|_0 < \text{spark}(\mathbf{X})/2$.

3.2. *$L_2$ penalty.* When the $\rho$ penalty is taken to be the $L_2$ penalty in the $\rho$-regularization problem (3), its minimizer is given by

$$\boldsymbol{\beta}_2 = \arg\min_{\boldsymbol{\beta} \in \mathcal{A}} \|\boldsymbol{\beta}\|_2. \tag{15}$$

It admits a closed-form solution. Viewing in the linear regression setting, we know that [see Theorem 6.2.1 in Fang and Zhang (1990)] the least squares estimate $(\mathbf{X}^T\mathbf{X})^+\mathbf{X}^T\mathbf{y}$ is a solution to the normal equation $\mathbf{X}^T\mathbf{y} = \mathbf{X}^T\mathbf{X}\boldsymbol{\beta}$, where $(\cdot)^+$ denotes the Moore–Penrose generalized matrix inverse. The following proposition shows that it coincides with $\boldsymbol{\beta}_2$.

PROPOSITION 1 (Minimum $L_2$ solution). $\boldsymbol{\beta}_2 = (\mathbf{X}^T\mathbf{X})^+\mathbf{X}^T\mathbf{y}$.

However, the minimum $L_2$ solution $\boldsymbol{\beta}_2$ to (2) is generally nonsparse and, thus, is different from the minimum $L_0$ solution $\boldsymbol{\beta}_0$.

3.3. *Penalties satisfying Condition 1.* We are curious about the $\rho/L_0$ equivalence, in the sense that the minimizer of the $\rho$-regularization problem (3) meets the minimum $L_0$ solution $\boldsymbol{\beta}_0$. As mentioned in the Introduction, many researchers have contributed to identifying conditions that ensure the $L_1/L_0$ equivalence when $\rho$ is taken to be the $L_1$ penalty. We consider penalties $\rho$ satisfying Condition 1. It is generally difficult to study the global minimizer analytically without convexity. As is common in the literature, we study the behavior of local minimizers.

Directly studying the local minimizer of the $\rho$-regularization problem (3) is generally difficult. We take the idea of constructing a solution to (3) by analyzing the solution of related $\rho$-regularized least squares problem (4) with regularization parameter $\lambda \in (0, \infty)$ and then letting $\lambda \to 0+$, where $p_\lambda(t) = \lambda \rho(t)$, $t \in [0, \infty)$.

We introduce some notation to simplify our presentation. For any $S \subset \{1, \ldots, p\}$, $\mathbf{X}_S$ stands for an $n \times |S|$ submatrix of $\mathbf{X}$ formed by columns with indices in $S$, $\mathbf{b}_S$ stands for the subvector of $\mathbf{b}$ formed by components with indices in $S$ and $S^c$ denotes its complement. For any vector $\mathbf{b} = (b_1, \ldots, b_q)^T$, define $\text{sgn}(\mathbf{b}) = (\text{sgn}(b_1), \ldots, \text{sgn}(b_q))^T$, where the sign function $\text{sgn}(x) = 1$ if $x > 0$, $-1$ if $x < 0$ and $0$ if $x = 0$. Let

$$\bar{\rho}(t) = \text{sgn}(t)\rho'(|t|), \qquad t \in \mathbf{R}, \tag{16}$$

and $\bar{\rho}(\mathbf{b}) = (\bar{\rho}(b_1), \ldots, \bar{\rho}(b_q))^T$, $\mathbf{b} = (b_1, \ldots, b_q)^T$. Clearly, for $a \in (0, \infty)$,

(17)  $\bar{\rho}_a(t) = \text{sgn}(t)a(a+1)/(a+|t|)^2$  and  $\bar{\rho}_\infty(t) = \text{sgn}(t), \qquad t \in \mathbf{R},$

where $\rho_a$ is defined in (9) and (10).

The following theorem gives a sufficient condition on the strict local minimizer of (4) for any $n$-vector $\mathbf{y}$ and $n \times p$ matrix $\mathbf{X}$.

THEOREM 1 (Regularized least squares). *Assume that $p_\lambda$ satisfies Condition 1 and $\widehat{\boldsymbol{\beta}}^\lambda \in \mathbf{R}^p$ with $\mathbf{Q} = \mathbf{X}_{\widehat{\mathfrak{M}}_\lambda}^T \mathbf{X}_{\widehat{\mathfrak{M}}_\lambda}$ nonsingular, where $\lambda \in (0, \infty)$ and $\widehat{\mathfrak{M}}_\lambda = \operatorname{supp}(\widehat{\boldsymbol{\beta}}^\lambda)$. Then $\widehat{\boldsymbol{\beta}}^\lambda$ is a strict local minimizer of (4) with $\lambda_n = \lambda$ if*

$$\widehat{\boldsymbol{\beta}}_{\widehat{\mathfrak{M}}_\lambda}^\lambda = \mathbf{Q}^{-1} \mathbf{X}_{\widehat{\mathfrak{M}}_\lambda}^T \mathbf{y} - \Lambda_n \lambda \mathbf{Q}^{-1} \bar{\rho}(\widehat{\boldsymbol{\beta}}_{\widehat{\mathfrak{M}}_\lambda}^\lambda), \tag{18}$$

$$\|\mathbf{z}_{\widehat{\mathfrak{M}}_\lambda^c}\|_\infty < \rho'(0+), \tag{19}$$

$$\lambda_{\min}(\mathbf{Q}) > \Lambda_n \lambda \kappa(\rho; \widehat{\boldsymbol{\beta}}_{\widehat{\mathfrak{M}}_\lambda}^\lambda), \tag{20}$$

*where $\mathbf{z} = (\Lambda_n \lambda)^{-1} \mathbf{X}^T (\mathbf{y} - \mathbf{X} \widehat{\boldsymbol{\beta}}^\lambda)$, $\lambda_{\min}(\cdot)$ denotes the smallest eigenvalue of a given symmetric matrix, and $\kappa(\rho; \widehat{\boldsymbol{\beta}}_{\widehat{\mathfrak{M}}_\lambda}^\lambda)$ is given by (8).*

Conditions (18) and (20) ensure that $\widehat{\boldsymbol{\beta}}^\lambda$ is a strict local minimizer of (4) when constrained on the $\|\widehat{\boldsymbol{\beta}}^\lambda\|_0$-dimensional subspace $\{\boldsymbol{\beta} \in \mathbf{R}^p : \boldsymbol{\beta}_{\widehat{\mathfrak{M}}_\lambda^c} = \mathbf{0}\}$ of $\mathbf{R}^p$. Condition (19) makes sure that the sparse vector $\widehat{\boldsymbol{\beta}}^\lambda$ is indeed a strict local minimizer of (4) on the whole space $\mathbf{R}^p$. When $\rho$ is convex, (19) and (20) can be, respectively, relaxed to no greater than and no less than under which $\widehat{\boldsymbol{\beta}}^\lambda$ is a minimizer of (4). Due to the possible nonconvexity of $\rho$, the technical analysis for proving local minimizer needs the strict inequalities in (19) and (20).

When $\rho$ is taken to be the $L_1$ penalty, the objective function in (4) is convex. Then the classical convex optimization theory applies to show that $\widehat{\boldsymbol{\beta}}^\lambda = (\widehat{\beta}_1^\lambda, \ldots, \widehat{\beta}_p^\lambda)$ is a global minimizer if and only if there exists a subgradient $\mathbf{z} \in \partial L_1(\widehat{\boldsymbol{\beta}}^\lambda)$, such that

$$\mathbf{X}^T \mathbf{X} \widehat{\boldsymbol{\beta}}^\lambda - \mathbf{X}^T \mathbf{y} + \Lambda_n \lambda \mathbf{z} = \mathbf{0}, \tag{21}$$

where the subdifferential of the $L_1$ penalty is given by $\partial L_1(\widehat{\boldsymbol{\beta}}^\lambda) = \{\mathbf{z} = (z_1, \ldots, z_p)^T \in \mathbf{R}^p : z_j = \operatorname{sgn}(\widehat{\beta}_j^\lambda) \text{ for } \widehat{\beta}_j^\lambda \neq 0 \text{ and } z_j \in [-1, 1] \text{ otherwise}\}$. Thus provided that $\mathbf{Q} = \mathbf{X}_{\widehat{\mathfrak{M}}_\lambda}^T \mathbf{X}_{\widehat{\mathfrak{M}}_\lambda}$ is nonsingular with $\widehat{\mathfrak{M}}_\lambda = \operatorname{supp}(\widehat{\boldsymbol{\beta}}^\lambda)$, condition (21) reduces to (18) and (19) with $<$ and $\rho'(0+)$ replaced by $\leq$ and $1$, respectively, whereas by the nonsingularity of $\mathbf{Q}$, condition (20) always holds for the $L_1$ penalty. However, to ensure that $\widehat{\boldsymbol{\beta}}^\lambda$ is the strict minimizer we need the strict inequality in (19). These conditions for the $L_1$ penalty have been extensively studied by many authors, for example, Efron et al. (2004), Fuchs (2004), Tropp (2006), Wainwright (2006) and Zhao and Yu (2006), among others.

By analyzing the solution of (4) characterized by Theorem 1 and letting $\lambda \to 0+$, we obtain the following theorem providing a sufficient condition which ensures that $\boldsymbol{\beta}_0$ is a local minimizer of (3).

THEOREM 2 (Sparse recovery). *Assume that $\rho$ satisfies Condition 1 with $\kappa(\rho) \in [0, \infty)$, $\mathbf{Q} = \mathbf{X}_{\mathfrak{M}_0}^T \mathbf{X}_{\mathfrak{M}_0}$ is nonsingular with $\mathfrak{M}_0 = \mathrm{supp}(\boldsymbol{\beta}_0)$, and $\mathbf{X} = (\mathbf{x}_1, \ldots, \mathbf{x}_p)$. Then $\boldsymbol{\beta}_0$ is a local minimizer of (3) if there exists some $\epsilon \in (0, \min_{j \in \mathfrak{M}_0} |\beta_{0,j}|)$ such that*

$$(22) \qquad \max_{j \in \mathfrak{M}_0^c} \max_{\mathbf{u} \in \mathcal{U}_\epsilon} |\langle \mathbf{x}_j, \mathbf{u} \rangle| < \rho'(0+),$$

*where $\mathcal{U}_\epsilon = \{\mathbf{X}_{\mathfrak{M}_0} \mathbf{Q}^{-1} \bar{\rho}(\mathbf{v}) : \mathbf{v} \in \mathcal{V}_\epsilon\}$ and $\mathcal{V}_\epsilon = \prod_{j \in \mathfrak{M}_0} \{t : |t - \beta_{0,j}| \le \epsilon\}$.*

We make some remarks on Theorem 2. Clearly condition (22) is free of the scale of the penalty function, that is, invariant under the rescaling $\rho \to c\rho$ for any constant $c \in (0, \infty)$. Condition (22) is independent of the scale of $\mathbf{X}$, and depends on $\boldsymbol{\beta}_0$ through $\mathfrak{M}_0$ and a small neighborhood $\mathcal{V}_\epsilon$ of $\boldsymbol{\beta}_{0,\mathfrak{M}_0}$. It can be viewed as a local condition at $(\mathfrak{M}_0, \boldsymbol{\beta}_{0,\mathfrak{M}_0})$. For the $L_1$ penalty $\rho_\infty$, we have $\rho'_\infty(0+) = 1$ and for any $\epsilon \in (0, \min_{j \in \mathfrak{M}_0} |\beta_{0,j}|)$, $\mathcal{U}_\epsilon$ contains a single point $\mathbf{u}_0 = \mathbf{X}_{\mathfrak{M}_0} \mathbf{Q}^{-1} \mathrm{sgn}(\boldsymbol{\beta}_{0,\mathfrak{M}_0})$. This shows that for the $L_1$ penalty, condition (22) reduces to the following condition

$$(23) \qquad \max_{j \in \mathfrak{M}_0^c} |\langle \mathbf{x}_j, \mathbf{u}_0 \rangle| < 1,$$

which can actually be relaxed to $\max_{j \in \mathfrak{M}_0^c} |\langle \mathbf{x}_j, \mathbf{u}_0 \rangle| \le 1$ while ensuring the $L_1/L_0$ equivalence. In the context of model selection, condition (23) was called the weak irrepresentable condition in Zhao and Yu (2006), who introduced it for characterizing the selection consistency of lasso. For the $\rho_a$ penalty with $a \in (0, \infty)$, by (11) we have

$$\rho'_a(0+) = 1 + a^{-1} \to \infty \quad \text{and} \quad \rho'_a(t) = \frac{a(a+1)}{(a+t)^2} \to 0 \qquad \text{as } a \to 0+$$

for each $t \in (0, \infty)$, which shows that condition (22) is less restrictive for smaller $a$. This justifies the flexibility of the $\rho_a$ penalties.

A great deal of research has contributed to identifying conditions on $\mathbf{X}$ and $\boldsymbol{\beta}_0$ that ensure the $L_1/L_0$ equivalence. See, for example, Chen, Donoho and Saunders (1999), Donoho and Elad (2003), Donoho (2004) and Candes and Tao (2005, 2006). In particular, Donoho (2004) shows that the individual equivalence of $L_1/L_0$ depends only on $\mathfrak{M}_0$ and $\boldsymbol{\beta}_{0,\mathfrak{M}_0}$. Condition (23) is independent of the scale of $\mathbf{X}$, and depends on $\boldsymbol{\beta}_0$ only through $\mathfrak{M}_0$ and $\boldsymbol{\beta}_{0,\mathfrak{M}_0}$, sharing the same spirit. The idea of using weighted $L_1$ penalty in the $\rho$-regularization problem (3) has been proposed and studied by Candès, Wakin and Boyd (2008).

3.4. *Optimal $\rho_a$ penalty.* Theorem 2 gives one characterization of the role of penalty functions in sparse recovery (2). In this section, we identify the optimal penalty $\rho_a$ for given $\mathbf{X}$ and $\boldsymbol{\beta}_0$ in sparse recovery.

For any $\epsilon \in (0, \min_{j \in \mathfrak{M}_0} |\beta_{0,j}|)$, we define

(24) $\quad \mathcal{P}_\epsilon = \{\text{All penalties } \rho_a \text{ in (3) that satisfy condition (22)}\}.$

By Theorem 2, any $\rho_a$ penalty in $\mathcal{P}_\epsilon$ ensures that $\boldsymbol{\beta}_0$ is recoverable in theory by the $\rho_a$-regularization problem (3). We are interested in a penalty $\rho_{a_{\mathrm{opt}}(\epsilon)}$ that attains the minimal maximum concavity in the sense that

(25) $$\kappa(\rho_{a_{\mathrm{opt}}(\epsilon)}) = \inf_{\rho_a \in \mathcal{P}_\epsilon} \kappa(\rho_a).$$

Such penalty $\rho_{a_{\mathrm{opt}}(\epsilon)}$ makes the objective function in (3) have the minimal maximum concavity, which is favorable from the computational point of view since the degree of concavity is related to the computational difficulty. We thus call $\rho_{a_{\mathrm{opt}}(\epsilon)}$ the optimal penalty. The following theorem characterizes it.

THEOREM 3 (Optimal $\rho_a$ penalty for sparse recovery). *Assume that $\mathbf{Q} = \mathbf{X}_{\mathfrak{M}_0}^T \mathbf{X}_{\mathfrak{M}_0}$ is nonsingular with $\mathfrak{M}_0 = \mathrm{supp}(\boldsymbol{\beta}_0)$ and $\epsilon \in (0, \min_{j \in \mathfrak{M}_0} |\beta_{0,j}|)$. Then the optimal penalty $\rho_{a_{\mathrm{opt}}(\epsilon)}$ satisfies:*

(a) $a_{\mathrm{opt}}(\epsilon) \in (0, \infty]$ *and is the largest $a \in (0, \infty]$ such that*

(26) $$\max_{j \in \mathfrak{M}_0^c} \max_{\mathbf{u} \in \mathcal{U}_\epsilon} |\langle \mathbf{x}_j, \mathbf{u} \rangle| \leq 1 + a^{-1},$$

*where $\mathcal{U}_\epsilon = \{\mathbf{X}_{\mathfrak{M}_0} \mathbf{Q}^{-1} \bar{\rho}(\mathbf{v}) : \mathbf{v} \in \mathcal{V}_\epsilon\}$ and $\mathcal{V}_\epsilon = \prod_{j \in \mathfrak{M}_0} \{t : |t - \beta_{0,j}| \leq \epsilon\}$.*
(b) $a_{\mathrm{opt}}(\epsilon) = \infty$ *if and only if*

(27) $$\max_{j \in \mathfrak{M}_0^c} |\langle \mathbf{x}_j, \mathbf{u}_0 \rangle| \leq 1,$$

*where $\mathbf{u}_0 = \mathbf{X}_{\mathfrak{M}_0} \mathbf{Q}^{-1} \mathrm{sgn}(\boldsymbol{\beta}_{0,\mathfrak{M}_0})$.*

By the characterization of $a_{\mathrm{opt}}(\epsilon)$, we see that for any $\rho_a$ penalty with $a \in (0, a_{\mathrm{opt}}(\epsilon))$, $\boldsymbol{\beta}_0$ is always a local minimizer of (3). However, in view of (12), its maximum concavity increases to $\infty$ as $a \to 0+$. Theorem 3(b) makes a sensible statement that for sparse recovery, the $L_1$ regularization in (3) is favorable from the computational point of view if condition (27) holds, which, as mentioned before, entails the $L_1/L_0$ equivalence. We would like to point out that the optimal parameter $a_{\mathrm{opt}}(\epsilon)$ depends on $\boldsymbol{\beta}_0$ and thus should be learned from the data. We give an example of calculating $a_{\mathrm{opt}}(\epsilon)$ below.

EXAMPLE 1. Assume that $\mathbf{X}_{\mathfrak{M}_0}$ is orthonormal with $\mathfrak{M}_0 = \{1, \ldots, s\}$ and $|\beta_{0,1}| = \cdots = |\beta_{0,s}|$. Thus by (2),

$$\mathbf{y} = \sum_{j=1}^{s} \beta_{0,j} \mathbf{x}_j = |\beta_{0,1}| \sum_{j=1}^{s} \mathrm{sgn}(\beta_{0,j}) \mathbf{x}_j. \tag{28}$$

Let $\mathcal{H}$ be the linear subspace of $\mathbf{R}^n$ spanned by $\mathbf{x}_1, \ldots, \mathbf{x}_s$ and $\mathcal{H}^\perp$ its orthogonal complement. Further assume that $p \geq s+1$, $\mathbf{x}_j \in \mathcal{H}^\perp$ for each $j \geq s+2$, $\|\mathbf{x}_{s+1}\|_2 = 1$, and

$$\mathbf{x}_{s+1} = rs^{-1/2} \sum_{j=1}^{s} \mathrm{sgn}(\beta_{0,j}) \mathbf{x}_j + \mathbf{h} \quad \text{with } r \in (-1, 1) \text{ and } \mathbf{h} \in \mathcal{H}^\perp. \tag{29}$$

Then for any $\epsilon \in (0, |\beta_{0,1}|)$, we have

$$\max_{j \in \mathfrak{M}_0^c} \max_{\mathbf{u} \in \mathcal{U}_\epsilon} |\langle \mathbf{x}_j, \mathbf{u} \rangle| = \max_{\mathbf{u} \in \mathcal{U}_\epsilon} |\langle \mathbf{x}_{s+1}, \mathbf{u} \rangle| = \max_{\mathbf{v} \in \mathcal{V}_\epsilon} |\langle \mathbf{x}_{s+1}, \mathbf{X}_{\mathfrak{M}_0} \bar{\rho}_a(\mathbf{v}) \rangle|$$

$$= |r| \sqrt{s} a(a+1) / (a + |\beta_{0,1}| - \epsilon)^2.$$

Thus condition (26) reduces to $|r|\sqrt{s} \frac{a(a+1)}{(a+|\beta_{0,1}|-\epsilon)^2} \leq 1 + a^{-1}$, which shows that for any $|r| > s^{-1/2}$, we have

$$a_{\mathrm{opt}}(\epsilon) = \frac{|\beta_{0,1}| - \epsilon}{(r^2 s)^{1/4} - 1}. \tag{30}$$

We see that the optimal parameter $a_{\mathrm{opt}}(\epsilon)$ is related to $|\beta_{0,1}| - \epsilon$ through both $r$ and $s$. It approaches $\infty$ as $|r| \to s^{-1/2}+$ and goes to $(|\beta_{0,1}| - \epsilon)/(s^{1/4} - 1)$ as $|r| \to 1-$ (see Figure 2). When $|r| \leq s^{-1/2}$, we have $a_{\mathrm{opt}}(\epsilon) = \infty$ regardless of $\epsilon \in (0, |\beta_{0,1}|)$.

In light of (28) and (29), $r \in (-1, 1)$ defining the noise predictor $\mathbf{x}_{s+1}$ is exactly the correlation between $\mathbf{x}_{s+1}$ and $\mathbf{y}$. Therefore in the noiseless setting (2), when the number of true variables $s$ is large, the correlation $r$ between the noise variable $\mathbf{x}_{s+1}$ and the response $\mathbf{y}$ has to be small in magnitude in order that the $L_1$ penalty is the optimal penalty, that is, $a_{\mathrm{opt}}(\epsilon) = \infty$. Note that the cut-off point for $|r|$ by our theory is $s^{-1/2}$, while $s^{-1/2}$ is exactly the absolute correlation between each true variable $\mathbf{x}_j$ and response $\mathbf{y}$.

If we assume $|\beta_{0,1}| = \min_{j \in \mathfrak{M}_0} |\beta_{0,j}|$ instead of $|\beta_{0,j}| = \cdots = |\beta_{0,s}|$, then the right-hand side of (30) gives a lower bound on $a_{\mathrm{opt}}(\epsilon)$, which entails that the cut-off point for $|r|$ by our theory can be above $s^{-1/2}$. This result is sensible once we observe that

$$|\mathrm{corr}(\mathbf{x}_{s+1}, \mathbf{y})| \leq |r| \quad \text{and} \quad \max_{1 \leq j \leq s} |\mathrm{corr}(\mathbf{x}_j, \mathbf{y})| = \frac{\max_{j \in \mathfrak{M}_0} |\beta_{0,j}|}{\sqrt{\sum_{j=1}^{s}(\beta_{0,j})^2}} \geq s^{-1/2},$$

where $\mathbf{y} = \sum_{i=1}^{s} \beta_{0,j} \mathbf{x}_j$. It is interesting to observe that the right-hand side of (30) is positively related to $\min_{j \in \mathfrak{M}_0} |\beta_{0,j}|$, which measures the strength of the weakest additive component in the response.

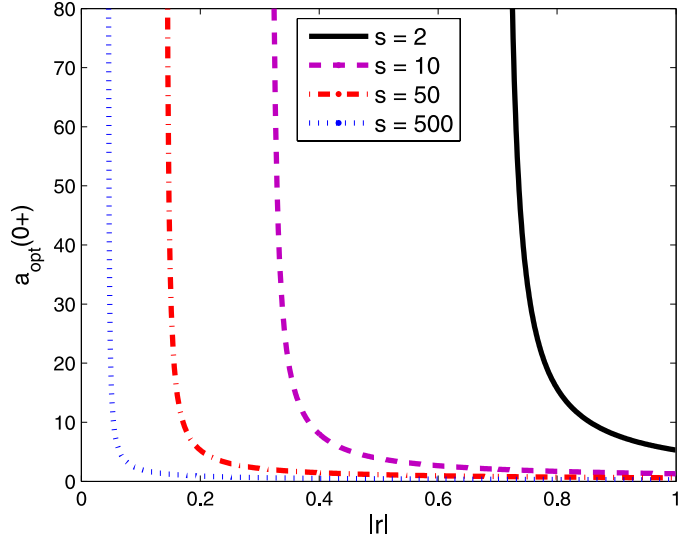

FIG. 2. *Plot of optimal parameter $a_{\mathrm{opt}}(0+)$ in (30) with $|\beta_{0,1}| = 1$ against $|r|$ for $s = 2$ (solid), $s = 10$ (dashed), $s = 50$ (dash-dot) and $s = 500$ (dotted).*

**4. Model selection.** In this section we consider the regularized least squares problem (4) for model selection in (1). Difficulties of recovering the true underlying sparse model $\mathfrak{M}_0$ include the collinearity among the predictors and the computational cost, both of which increase as the dimensionality grows [see, e.g., Fan and Li (2006) and Fan and Lv (2008)]. For example, classical model selection approaches such as best subset selection are very demanding in computation and become impractical to implement in high dimensions.

There is a huge literature on model selection. To name a few in addition to those work mentioned before, Frank and Friedman (1993) propose the bridge regression. Breiman (1995) introduces the nonnegative garrote for shrinkage estimation and variable selection. Tibshirani (1996) proposes the lasso using the $L_1$-regularized least squares. Oracle properties of nonconcave penalized likelihood estimators including the SCAD [Fan and Li (2001)] have been systematically studied by Fan and Li (2001) and Fan and Peng (2004). In particular, Fan and Li (2001) propose a unified algorithm LQA for optimizing nonconcave penalized likelihood. Efron et al. (2004) introduce the least angle regression for variable selection and present the LARS algorithm. Zou and Li (2008) propose one-step sparse estimates for nonconcave penalized likelihood models and introduce the LLA algorithm for optimizing nonconcave penalized likelihood. Candes and Tao (2007) propose the Dantzig selector and prove its nonasymptotic properties. Later, Meinshausen, Rocha and Yu (2007), Bickel, Ritov and Tsybakov (2008) and James, Radchenki and Lv (2009) establish the equivalence or approximate equivalence of lasso and Dantzig selector under different conditions. More recent regularization methods include MCP proposed by Zhang (2007).

Fan and Li (2001) point out the bias issue of lasso. Zou (2006) proposes the adaptive lasso to address this issue by using an adaptively weighted $L_1$ penalty. Greenshtein and Ritov (2004) study the persistency of lasso-type procedures in high-dimensional linear predictor selection. Hunter and Li (2005) propose and study MM algorithms for variable selection. Li and Liang (2008) study variable selection for semiparametric regression models. Wang, Li and Tsai (2007) study the problem of tuning parameter selection for the SCAD. Fan and Fan (2008) study the impact of high dimensionality on classifications and propose the FAIR. Fan and Lv (2008) propose the SIS as well as its extensions for variable screening and study its asymptotic properties in ultra-high-dimensional feature space.

4.1. *Regularized least squares estimator.* We consider the regularized least squares problem (4) with the penalty $p_\lambda$ in the class satisfying Condition 1. For a given regularization parameter $\lambda_n \in [0, \infty)$ indexed by sample size $n$, a $p$-dimensional vector $\widehat{\boldsymbol{\beta}}^{\lambda_n}$ is conventionally called a regularized least squares estimator of $\boldsymbol{\beta}_0$ if it is a (local) minimizer of (4). When the $L_2$ penalty $\rho(t) = t^2$ is used, the resulting estimator is called the ridge estimator, and its limit as $\lambda_n \to 0+$ can be easily shown to be the ordinary least squares estimator $\widehat{\boldsymbol{\beta}}^{\text{ols}} \equiv (\mathbf{X}^T\mathbf{X})^+\mathbf{X}^T\mathbf{y}$, where $(\cdot)^+$ denotes the Moore–Penrose generalized matrix inverse. $\widehat{\boldsymbol{\beta}}^{\text{ols}}$ is also a solution to the normal equation $\mathbf{X}^T\mathbf{y} = \mathbf{X}^T\mathbf{X}\boldsymbol{\beta}$.

When $\lambda_n \in (0, \infty)$, Theorem 1 in Section 3.3 gives a sufficient condition on the strict local minimizer of (4). From the proof of Theorem 1, we see that any local minimizer $\widehat{\boldsymbol{\beta}}^{\lambda_n}$ of (4) must satisfy

$$\widehat{\boldsymbol{\beta}}^{\lambda_n}_{\widehat{\mathfrak{M}}_{\lambda_n}} = \mathbf{Q}^{-1}\mathbf{X}^T_{\widehat{\mathfrak{M}}_{\lambda_n}}\mathbf{y} - \Lambda_n\lambda_n\mathbf{Q}^{-1}\bar{\rho}(\widehat{\boldsymbol{\beta}}^{\lambda_n}_{\widehat{\mathfrak{M}}_{\lambda_n}}), \tag{31}$$

$$\|\mathbf{z}_{\widehat{\mathfrak{M}}^c_{\lambda_n}}\|_\infty \leq \rho'(0+), \tag{32}$$

$$\lambda_{\min}(\mathbf{Q}) \geq \Lambda_n\lambda_n\kappa(\rho; \widehat{\boldsymbol{\beta}}^{\lambda_n}_{\widehat{\mathfrak{M}}_{\lambda_n}}), \tag{33}$$

where $\widehat{\mathfrak{M}}_{\lambda_n} = \text{supp}(\widehat{\boldsymbol{\beta}}^{\lambda_n})$, $\mathbf{Q} = \mathbf{X}^T_{\widehat{\mathfrak{M}}_{\lambda_n}}\mathbf{X}_{\widehat{\mathfrak{M}}_{\lambda_n}}$, $\mathbf{z} = (\Lambda_n\lambda_n)^{-1}\mathbf{X}^T(\mathbf{y} - \mathbf{X}\widehat{\boldsymbol{\beta}}^{\lambda_n})$, and $\kappa(\rho; \widehat{\boldsymbol{\beta}}^{\lambda_n}_{\widehat{\mathfrak{M}}_{\lambda_n}})$ is given by (8). So there is generally a slight gap between the necessary condition for local minimizer and sufficient condition for strict local minimizer, in view of (32), (33) and (19), (20). Hereafter the regularized least squares estimator is referred to as a $Z$-estimator $\widehat{\boldsymbol{\beta}}^{\lambda_n} \in \mathbf{R}^p$ that solves (31)–(33).

We observe that $\widehat{\boldsymbol{\beta}}^{\lambda_n}_{\widehat{\mathfrak{M}}_{\lambda_n}}$ in (31) is the difference between two terms. The first term $\mathbf{Q}^{-1}\mathbf{X}^T_{\widehat{\mathfrak{M}}_{\lambda_n}}\mathbf{y}$ is exactly the ordinary least squares estimator by using predictors $\mathbf{x}_j$ with indices in $\widehat{\mathfrak{M}}_{\lambda_n}$. In the case of orthonormal design matrix $\mathbf{X}$, we have $\mathbf{Q} = I_{\widehat{s}_n}$ with $\widehat{s}_n = \|\widehat{\boldsymbol{\beta}}^{\lambda_n}\|_0$, and thus the second term becomes $\Lambda_n\lambda_n\bar{\rho}(\widehat{\boldsymbol{\beta}}^{\lambda_n}_{\widehat{\mathfrak{M}}_{\lambda_n}})$, which, for nonzero components, has the same sign as $\widehat{\boldsymbol{\beta}}^{\lambda_n}_{\widehat{\mathfrak{M}}_{\lambda_n}}$ componentwise by definition. In view of (31), we have

$$\widehat{\boldsymbol{\beta}}^{\lambda_n}_{\widehat{\mathfrak{M}}_{\lambda_n}} + \Lambda_n\lambda_n\bar{\rho}(\widehat{\boldsymbol{\beta}}^{\lambda_n}_{\widehat{\mathfrak{M}}_{\lambda_n}}) = \mathbf{X}^T_{\widehat{\mathfrak{M}}_{\lambda_n}}\mathbf{y},$$

which entails that both $\widehat{\boldsymbol{\beta}}^{\lambda_n}_{\widehat{\mathfrak{M}}_{\lambda_n}}$ and $\Lambda_n\lambda_n\bar{\rho}(\widehat{\boldsymbol{\beta}}^{\lambda_n}_{\widehat{\mathfrak{M}}_{\lambda_n}})$ (for its nonzero components) have the same sign as the ordinary least squares estimator $\mathbf{X}^T_{\widehat{\mathfrak{M}}_{\lambda_n}}\mathbf{y}$. This shows that the second term above is indeed a shrinkage towards zero when $\mathbf{X}$ is orthonormal. For the penalties $\rho_a$ introduced in Section 2.2, $\bar{\rho}_a(t)$ depends on both $t$ and $a$ [see (17)]. In fact, for small $a$, $\bar{\rho}_a(t)$ takes a wide range of values when $t$ varies, which ensures that $\rho_a$ penalty shrinks the components of the ordinary least squares estimator differently. This gives us flexibility in model selection. It provides us a family of regularized least squares estimators indexed by parameter $a$ and regularization parameter $\lambda_n$. As $a$ becomes smaller, it generally gives sparser estimates.

4.2. *Nonasymptotic properties.* In this paper, we study a nonasymptotic property of $\widehat{\boldsymbol{\beta}}^{\lambda_n}$, called the weak oracle property for simplicity, which means: (1) sparsity in the sense of $\widehat{\boldsymbol{\beta}}^{\lambda_n}_{\mathfrak{M}_0^c} = \mathbf{0}$ with probability tending to 1 as $n \to \infty$, and (2) consistency under the $L_\infty$ loss. This property is weaker than the oracle property introduced by Fan and Li (2001). As mentioned before, we condition on the $n \times p$ design matrix $\mathbf{X}$.

We use the $p_\lambda$ penalty in the class satisfying Condition 1 and make the following assumptions on the deterministic design matrix $\mathbf{X}$ and the distribution of the noise vector $\boldsymbol{\varepsilon}$ in the linear model (1).

CONDITION 2. $\mathbf{X}$ satisfies

(34) $$\|(\mathbf{X}^T_{\mathfrak{M}_0}\mathbf{X}_{\mathfrak{M}_0})^{-1}\|_\infty \leq C_{1n},$$

(35) $$\|\mathbf{X}^T_{\mathfrak{M}_0^c}\mathbf{X}_{\mathfrak{M}_0}(\mathbf{X}^T_{\mathfrak{M}_0}\mathbf{X}_{\mathfrak{M}_0})^{-1}\|_\infty \leq C_{2n},$$

where $\mathfrak{M}_0 = \mathrm{supp}(\boldsymbol{\beta}_0)$, $C_{1n} \in (0, \infty)$, $C_{2n} \in [0, C\frac{\rho'(0+)}{\rho'(c_0 b_0)}]$ for some $C, c_0 \in (0, 1)$, $b_0 = \min_{j \in \mathfrak{M}_0}|\beta_{0,j}|$, and $\|\cdot\|_\infty$ denotes the matrix $\infty$-norm.

Here and below, $\rho$ is associated with regularization parameter $\underline{\lambda}_n$ defined in (38) unless specified otherwise.

CONDITION 3.  $\varepsilon \sim N(\mathbf{0}, \sigma^2 I_n)$ for some $\sigma > 0$.

When $\mathbf{X}_{\mathfrak{M}_0}$ is orthonormal, the left-hand side of (35) becomes $\|\mathbf{X}_{\mathfrak{M}_0^c}^T \mathbf{X}_{\mathfrak{M}_0}\|_\infty$, the maximum absolute sum of covariances between a noise variable and all $s$ true variables. Condition (35) constrains its growth rate. By the concavity of $\rho$ in Condition 1, $\rho'(t)$ is decreasing in $t \in [0, \infty)$ and thus the ratio $\rho'(0+)/\rho'(c_0 b_0)$ is always no less than one. When the $L_1$ penalty is used, $C_{2n} \in [0, C]$ and Condition 2 reduces to the condition in Zhao and Yu (2006) and Wainwright (2006).

Condition 2 is an assumption on design matrix $\mathbf{X}$. If we work with random design, we can calculate the probability that Condition 2 holds by using the results from, for example, the random matrix theory. The Gaussian assumption in Condition 3 can be relaxed to other light-tailed distributions so that we can derive similar exponential probability bounds.

CONDITION 4.  There exists some $\gamma \in (0, \frac{1}{2}]$ such that

$$(36) \qquad \left[D_{1n} + \frac{\rho'(c_0 b_0)}{\rho'(0+)} D_{2n}\right] C_{1n} = O(n^{-\gamma}),$$

where $D_{1n} = \max_{j \in \mathfrak{M}_0} \|\mathbf{x}_j\|_2$, $D_{2n} = \max_{j \in \mathfrak{M}_0^c} \|\mathbf{x}_j\|_2$ and $\mathbf{X} = (\mathbf{x}_1, \ldots, \mathbf{x}_p)$. Let $u_n \in (0, \infty)$ satisfy $\lim_{n \to \infty} u_n = \infty$, $\underline{\lambda}_n \leq \overline{\lambda}_n$, and

$$(37) \quad u_n \leq [\kappa_0 (C_{2n} D_{1n} + D_{2n})]^{-1} \lambda_{\min}(\mathbf{X}_{\mathfrak{M}_0}^T \mathbf{X}_{\mathfrak{M}_0})(1-C)\rho'(0+)\sigma^{-1},$$

where

$$(38) \quad \underline{\lambda}_n = \Lambda_n^{-1} \frac{(C_{2n} D_{1n} + D_{2n}) u_n \sigma}{\rho'(0+) - C_{2n} \rho'(c_0 b_0)} \quad \text{and} \quad \overline{\lambda}_n = \frac{C_{1n}^{-1}(1-c_0)b_0 - u_n D_{1n} \sigma}{\Lambda_n \rho'(c_0 b_0; \overline{\lambda}_n)},$$

$C, c_0 \in (0, 1)$ are given in Condition 2, and $\kappa_0 = \max\{\kappa(\rho; \mathbf{b}) : \|\mathbf{b} - \boldsymbol{\beta}_{0, \mathfrak{M}_0}\|_\infty \leq (1 - c_0) b_0\}$ with $\kappa(\rho; \mathbf{b})$ given by (8).

As seen later, we need the condition $\underline{\lambda}_n \leq \overline{\lambda}_n$ to ensure the existence of a desired regularization parameter $\lambda_n = \underline{\lambda}_n$. Condition (37), which always holds when $\kappa_0 = 0$, is needed to ensure condition (33).

THEOREM 4 (Weak oracle property).  *Assume that $p_\lambda$ in (4) satisfies Condition 1, Conditions 2–4 hold and $p = o(u_n e^{u_n^2/2})$. Then there exists a regularized least squares estimator $\widehat{\boldsymbol{\beta}}^{\lambda_n}$ with regularization parameter $\lambda_n = \underline{\lambda}_n$ defined in (38) such that with probability at least $1 - \frac{2}{\sqrt{\pi}} p u_n^{-1} e^{-u_n^2/2}$, $\widehat{\boldsymbol{\beta}}^{\lambda_n}$ satisfies:*

(a) (Sparsity) $\widehat{\boldsymbol{\beta}}_{\mathfrak{M}_0^c}^{\lambda_n} = \mathbf{0}$;

(b) ($L_\infty$ loss) $\|\widehat{\boldsymbol{\beta}}_{\mathfrak{M}_0}^{\lambda_n} - \boldsymbol{\beta}_{0,\mathfrak{M}_0}\|_\infty \leq h = O(n^{-\gamma} u_n)$,

where $\mathfrak{M}_0 = \mathrm{supp}(\boldsymbol{\beta}_0)$ and $h = [D_{1n} + \frac{\rho'(c_0 b_0)}{\rho'(0+)} D_{2n}] C_{1n} u_n (1-C)^{-1} \sigma$. As a consequence, $\|\widehat{\boldsymbol{\beta}}^{\lambda_n} - \boldsymbol{\beta}_0\|_2 = O_P(\sqrt{s} n^{-\gamma} u_n)$, where $s = \|\boldsymbol{\beta}_0\|_0$.

We make some remarks on Theorem 4. In view of $\lim_{n \to \infty} u_n = \infty$ in Condition 4, the dimensionality $p$ is allowed to grow up exponentially fast with $u_n$. If $\kappa_0$ in (37) is of a small order, $u_n$ can be allowed to be $o(n^\gamma)$ and thus $\log p = o(n^{2\gamma})$. The diverging sequence $(u_n)$ also controls the rate of the exponential probability bound. From Theorem 4, we see that with asymptotic probability one the $L_\infty$ estimation loss of $\widehat{\boldsymbol{\beta}}^{\lambda_n}$ is bounded above by $h_1 + h_2$, where

$$(39) \quad h_1 = D_{1n} C_{1n} u_n (1-C)^{-1} \sigma \quad \text{and} \quad h_2 = \frac{\rho'(c_0 b_0)}{\rho'(0+)} D_{2n} C_{1n} u_n (1-C)^{-1} \sigma.$$

The second term $h_2$ is associated with the penalty function $\rho$. For the $L_1$ penalty, the ratio $\rho'(c_0 b_0)/\rho'(0+)$ is equal to one, and for other concave penalties, as mentioned before, this ratio can be (much) smaller than one. This is consistent with the fact shown by Fan and Li (2001) that concave penalties can reduce the biases of estimates.

In the classical setting of $D_{1n}, D_{2n} = O(\sqrt{n})$ and $C_{1n} = O(n^{-1})$, we have $\gamma = 1/2$ since $\rho'(c_0 b_0)/\rho'(0+) \leq 1$ by Condition 1 and thus the consistency rate of $\widehat{\boldsymbol{\beta}}^{\lambda_n}$ under the $L_2$ norm becomes $O_P(\sqrt{s} n^{-1/2} u_n)$, which is slightly slower than $O_P(\sqrt{s} n^{-1/2})$. This is because it is derived by using the $L_\infty$ loss of $\widehat{\boldsymbol{\beta}}^{\lambda_n}$ in Theorem 4(b). The use of the $L_\infty$ norm is due to the technical difficulty of proving the existence of a solution to the nonlinear equation (31). We conjecture that by considering the $L_2$ norm directly, one can obtain the consistency rate $O_P(\sqrt{s} n^{-1/2})$.

4.3. *Choice of $\rho_a$ penalty.* Theorem 4 gives one characterization of the role of penalty functions in regularized least squares for model selection in (1). Let us now consider the penalties $\rho_a$ introduced in Section 2.2. We fix the diverging sequence $(u_n)$ in Theorem 4 and see how the parameter $a \in (0, \infty]$ influences the performance of the $\rho_a$-regularized least squares method.

In view of (11), we have

$$(40) \quad \frac{\rho_a'(0+)}{\rho_a'(c_0 b_0)} = \frac{(a + c_0 b_0)^2}{a^2}, \quad a \in (0, \infty) \quad \text{and} \quad \frac{\rho_\infty'(0+)}{\rho_\infty'(c_0 b_0)} = 1.$$

Clearly $\rho_a'(0+)/\rho_a'(c_0 b_0)$ is decreasing in $a \in (0, \infty]$. Thus (35) in Condition 2 becomes less restrictive as $a \to 0+$. In view of (39) by Theorem 4, the upper bound on the $L_\infty$ estimation loss of the $\rho_a$-regularized least squares estimator $\widehat{\boldsymbol{\beta}}^{\lambda_n}$ decreases to $h_1$ as $a$ approaches 0. However, in view of (12), the maximum concavity of $\rho_a$ increases to $\infty$ as $a \to 0+$. This suggests that the computational difficulty of solving the $\rho_a$-regularized least squares problem (4) may increase as $a$ approaches 0. In practical implementation, we can adaptively choose $a$ using the data, for example, the cross-validation method.

**5. Implementation.** In this section, we discuss algorithms for solving regularization problems with $\rho_a$ penalty. Specifically, the $\rho_a$-regularization problem (3) for sparse recovery in (2), and the $\rho_a$-regularized least squares problem (4) for model selection in (1).

5.1. *SIRS algorithm for sparse recovery.* For any $\boldsymbol{\beta} = (\beta_1, \ldots, \beta_p)^T \in \mathbf{R}^p$, let $\mathbf{D}(\boldsymbol{\beta}) = \mathrm{diag}\{d_1, \ldots, d_p\}$ and

$$\mathbf{v}(\boldsymbol{\beta}) = \mathbf{D}\mathbf{X}^T(\mathbf{X}\mathbf{D}\mathbf{X}^T)^+\mathbf{y}, \tag{41}$$

where $\mathbf{D} = \mathbf{D}(\boldsymbol{\beta})$ and $d_j = \beta_j^2/\rho_a(|\beta_j|) = (a+1)^{-1}|\beta_j|(a+|\beta_j|)$, $j = 1, \ldots, p$. We propose the sequentially and iteratively reweighted squares (SIRS) algorithm for solving (3) with $\rho_a$ penalty. SIRS uses the method of iteratively reweighted squares and iteratively solves the $\rho$-regularization problem (3) with the weighted $L_2$ penalty $\rho(\boldsymbol{\beta}) = \boldsymbol{\beta}^T\boldsymbol{\Gamma}\boldsymbol{\beta}$, where $\boldsymbol{\Gamma} = \mathrm{diag}\{d_1^{-1}, \ldots, d_p^{-1}\}$ with $\mathbf{D} = \mathbf{D}(\boldsymbol{\gamma})$ for some $\boldsymbol{\gamma}$, $0^{-1} = \infty$ and $0 \cdot \infty = 0$. It sequentially searches for a good initial value that leads to $\boldsymbol{\beta}_0$. Pick a level of sparsity $S$, the number of iterations $L$, the number of sequential steps $M \leq S$ and a small constant $\epsilon \in (0, 1)$.

SIRS ALGORITHM.

1. Set $k = 0$.
2. Initialize $\boldsymbol{\beta}^{(0)} = \mathbf{1}$ and set $\ell = 1$.
3. Set $\boldsymbol{\beta}^{(\ell)} \leftarrow \mathbf{v}(\boldsymbol{\beta}^{(\ell-1)})$ with $\mathbf{D} = \mathbf{D}(\boldsymbol{\beta}^{(\ell-1)})$ and $\ell \leftarrow \ell + 1$.
4. Repeat step 3 until convergence or $\ell = L + 1$. Denote by $\widetilde{\boldsymbol{\beta}}$ the resulting $p$-vector.
5. If $\|\widetilde{\boldsymbol{\beta}}\|_0 \leq S$, stop and return $\widetilde{\boldsymbol{\beta}}$. Otherwise, set $k \leftarrow k + 1$ and repeat steps 2–4 with $\boldsymbol{\beta}^{(0)} = I(|\widetilde{\boldsymbol{\beta}}| \geq \gamma_k) + \epsilon I(|\widetilde{\boldsymbol{\beta}}| < \gamma_k)$ and $\gamma_k$ the $k$th largest component of $|\widetilde{\boldsymbol{\beta}}|$, until stop or $k = M$. Return $\widetilde{\boldsymbol{\beta}}$.

In practice, we can set a small tolerance level for convergence and apply a hard thresholding with a sufficiently small threshold to $\widetilde{\boldsymbol{\beta}}$ to generate sparsity in step 4. The small constant $\epsilon \in (0, 1)$ is introduced to leverage the scoring of variables with indices in $\mathfrak{M}_0 = \mathrm{supp}(\boldsymbol{\beta}_0)$ by suppressing the noise variables. Through numerical implementations, we have found that SIRS is robust to the choice of $\epsilon$. SIRS algorithm stops once it finds a sufficiently sparse solution to (2). We use a grid search method to select the optimal parameter $a$ of penalty $\rho_a$ that produces the sparsest solution.

5.1.1. *Justification.* It is nontrivial to derive the convergence of SIRS algorithm analytically. In all our numerical implementations, we have found that it always converges. In this section, we give a partial justification of SIRS.

For each given $\boldsymbol{\beta}^{(\ell-1)}$, step 3 of SIRS solves the $\rho$-regularization problem (3) with the weighted $L_2$ penalty $\rho(\boldsymbol{\beta}) = \boldsymbol{\beta}^T \boldsymbol{\Gamma} \boldsymbol{\beta}$, where $\boldsymbol{\Gamma}$ is given by $\mathbf{D} = \mathbf{D}(\boldsymbol{\beta}^{(\ell-1)})$. This can be easily shown by Proposition 1 and the identity $\mathbf{v}(\boldsymbol{\beta}) = \mathbf{D}^{1/2}(\mathbf{D}^{1/2}\mathbf{X}^T\mathbf{X}\mathbf{D}^{1/2})^{+}\mathbf{D}^{1/2}\mathbf{X}^T\mathbf{y}$. It follows from the definitions of $\mathbf{D}$ and $\boldsymbol{\Gamma}$ that $\rho_a(\boldsymbol{\beta}^{(\ell-1)}) = (\boldsymbol{\beta}^{(\ell-1)})^T \boldsymbol{\Gamma} \boldsymbol{\beta}^{(\ell-1)} = \rho(\boldsymbol{\beta}^{(\ell-1)})$. Thus the two penalties agree at $\boldsymbol{\beta} = \boldsymbol{\beta}^{(\ell-1)}$. Quadratic approximations have been used in many iterative algorithms such as LQA in Fan and Li (2001).

The following proposition characterizes $\lim_{\ell \to \infty} \boldsymbol{\beta}^{(\ell)}$ when it exists.

PROPOSITION 2. (a) *If* $\lim_{\ell \to \infty} \boldsymbol{\beta}^{(\ell)}$ *exists, then it is a fixed point of the functional* $\mathcal{F} : \mathcal{A} \to \mathcal{A}$,

$$(42) \qquad \mathcal{F}(\boldsymbol{\gamma}) = \arg\min_{\boldsymbol{\beta} \in \mathcal{A}} \boldsymbol{\beta}^T \boldsymbol{\Gamma}(\boldsymbol{\gamma}) \boldsymbol{\beta},$$

*where* $\mathcal{A} = \{\boldsymbol{\beta} \in \mathbf{R}^p : \mathbf{y} = \mathbf{X}\boldsymbol{\beta}\}$ *and* $\boldsymbol{\Gamma}(\boldsymbol{\gamma})$ *denotes* $\boldsymbol{\Gamma}$ *given by* $\mathbf{D} = \mathbf{D}(\boldsymbol{\gamma})$.

(b) $\boldsymbol{\beta}_0$ *is always a fixed point of the functional* $\mathcal{F}$.

(c) *Assume that* $p > n$, $\mathrm{spark}(\mathbf{X}) = n + 1$ *and* $\|\boldsymbol{\beta}_0\|_0 < (n+1)/2$. *Then for any fixed point* $\boldsymbol{\beta}$ *of the functional* $\mathcal{F}$, *we have* $\boldsymbol{\beta} = \boldsymbol{\beta}_0$ *or* $\|\boldsymbol{\beta}\|_0 > (n+1)/2$.

5.1.2. *Computational complexity.* We now discuss the computational complexity of the SIRS algorithm. We first consider the case of $p \geq n$. Note that $\mathbf{D}\mathbf{X}^T(\mathbf{X}\mathbf{D}\mathbf{X}^T)^{+}\mathbf{y} = \lim_{\lambda \to 0+} \mathbf{D}\mathbf{X}^T(\lambda I_n + \mathbf{X}\mathbf{D}\mathbf{X}^T)^{-1}\mathbf{y}$. Thus in step 3, we can avoid calculating the generalized matrix inverse and only need to calculate the $p$-vector $\mathbf{D}\mathbf{X}^T(\lambda I_n + \mathbf{X}\mathbf{D}\mathbf{X}^T)^{-1}\mathbf{y}$ for a fixed sufficiently small $\lambda \in (0, \infty)$. Since the $p \times p$ matrix $\mathbf{D}$ is diagonal, the computational complexity of $\lambda I_n + \mathbf{X}\mathbf{D}\mathbf{X}^T$ is $O(n^2 p)$. It is known that inverting an $n \times n$ matrix is of computational complexity $O(n^3)$. This shows that the computational complexity of $(\lambda I_n + \mathbf{X}\mathbf{D}\mathbf{X}^T)^{-1}$ is $O(n^2 p)$ since $p \geq n$. Then it is easy to see that step 3 has computational complexity $O(n^2 p)$. Similarly, when $p \leq n$ we can derive that step 3 has computational complexity $O(np^2)$ by using the identity $\mathbf{v}(\boldsymbol{\beta}) = \mathbf{D}^{1/2}(\mathbf{D}^{1/2}\mathbf{X}^T\mathbf{X}\mathbf{D}^{1/2})^{+}\mathbf{D}^{1/2}\mathbf{X}^T\mathbf{y}$, which equals $\lim_{\lambda \to 0+} \mathbf{D}^{1/2}(\lambda I_p + \mathbf{D}^{1/2}\mathbf{X}^T\mathbf{X}\mathbf{D}^{1/2})^{-1}\mathbf{D}^{1/2}\mathbf{X}^T\mathbf{y}$. Therefore, the computational complexity of the main step, step 3, of the SIRS algorithm is $O(np(n \wedge p))$, which is the same as that of the LARS algorithm [Efron et al. (2004)] for model selection. For given number of iterations $L$ and number of sequential steps $M$, SIRS has computational complexity $O(np(n \wedge p)LM)$. We would like to point out that SIRS stops at any sequential step once it finds a sufficiently sparse solution to (2).

5.2. *Model selection.* Efficient algorithms for solving the regularized least squares problem (4) include the LQA proposed by Fan and Li (2001) and LLA introduced by Zou and Li (2008). In this paper we use LLA to solve (4) with $\rho_a$ penalty. For a fixed regularization parameter $\lambda_n$, LLA iteratively solves (4) by using local linear approximations of $\sum_{j=1}^{p} \rho_a(|\beta_j|)$. At a given $\boldsymbol{\beta}^{(0)} = (\beta_1^{(0)}, \ldots, \beta_p^{(0)})^T$, LLA approximates $\sum_{j=1}^{p} \rho_a(|\beta_j|)$ as

$$\sum_{j=1}^{p} [\rho_a(|\beta_j^{(0)}|) + \rho_a'(|\beta_j^{(0)}|)(|\beta_j| - |\beta_j^{(0)}|)].$$

Then LLA solves the weighted lasso ($L_1$-regularized least squares)

$$(43) \qquad \min_{\boldsymbol{\beta} \in \mathbf{R}^p} \left\{ 2^{-1} \|\mathbf{y} - \mathbf{X}\boldsymbol{\beta}\|_2^2 + \Lambda_n \lambda_n \sum_{j=1}^{p} w_j |\beta_j| \right\},$$

where $w_j = \rho_a'(|\beta_j^{(0)}|) = a(a+1)/(a + |\beta_j^{(0)}|)^2$, $j = 1, \ldots, p$. We use a grid search method to tune the parameter $a$ of penalty $\rho_a$.

## 6. Numerical examples.

6.1. *Simulation* 1. In this example, we demonstrate the performance of $\rho_a$ penalty in sparse recovery. We simulated 100 data sets from (2) with $(s, n, p) = (7, 35, 1000)$, $\mathfrak{M}_0 = \{1, 2, \ldots, 7\}$ and $\boldsymbol{\beta}_{0, \mathfrak{M}_0} = (1, -0.5, 0.7, -1.2, -0.9, 0.3, 0.55)^T$, for each of three levels of correlation $r$. Let $\Gamma_r$ be a $p \times p$ matrix with diagonal elements being 1 and off-diagonal elements being $r$. We chose $r = 0$, 0.2 and 0.5. The rows of $\mathbf{X}$ were first sampled as i.i.d. copies from $N(\mathbf{0}, \Gamma_r)$ and then each of its columns was rescaled to have unit $L_2$ norm.

As discussed in Section 5.1, we implemented $\rho_a$ regularization (3) using SIRS algorithm. We set $S = \lceil \frac{n}{2} \rceil$, $\epsilon = p^{-1}$ and $M = S$. The parameter $a$ was chosen to be in $\{0, 0.05, 0.1, 0.2, 0.3, 0.4, 0.6, 1, 2, 5\}$. For a comparison, we also implemented the $L_1$ regularization (3), which can easily be recast as a linear program. The optimal SICA refers to the method that tunes parameter $a$ of $\rho_a$ penalty by selecting the one that generates the sparsest solution. Since sparse recovery can be formulated as model selection, we also implemented $\rho_a$-regularized least squares (4) using LLA algorithm with parameters $\lambda$ and $a$ tuned by cross-validation, and compared it with SIRS.

Table 1 shows the comparison results. We see that $\rho_a$ with finite $a$ and optimal SICA significantly outperformed $L_1$. As $a$ gets larger, the performance of $\rho_a$ approaches that of $L_1$. When $a$ approaches zero, the success percentage first increases and then decreases. This suggests that the computational difficulty increases for very small $a$.

6.2. *Simulation* 2. In this example as well as the next two ones, we demonstrate the performance of $\rho_a$ penalty in model selection. The data were generated from model (1). We set $(n, p) = (100, 50)$ and chose the true regression coefficients vector $\boldsymbol{\beta}_0$ by setting $\mathfrak{M}_0 = \{1, 2, \ldots, 7\}$ and $\boldsymbol{\beta}_{0,\mathfrak{M}_0} = (1, -0.5, 0.7, -1.2, -0.9, 0.3, 0.55)^T$. The number of simulations was 100. For each simulated data set, the rows of $\mathbf{X}$ were sampled as i.i.d. copies from $N(\mathbf{0}, \Sigma_0)$ with $\Sigma_0 = (0.5^{|i-j|})_{i,j=1,\ldots,p}$, and $\boldsymbol{\varepsilon}$ was generated independently from $N(\mathbf{0}, \sigma^2 I_n)$. Two noise levels $\sigma = 0.3$ and $0.5$ were considered. We compared SICA with lasso, SCAD and MCP. Lasso was implement by LARS algorithm, and SCAD, MCP and SICA were implemented by LLA algorithm. The regularization parameters $\lambda$ and $a$ were selected by using a grid search method based on BIC, following Wang, Li and Tsai (2007).

Three performance measures were employed to compare the four methods. The first measure is the prediction error (PE) defined as $E(y - \mathbf{x}^T \widehat{\boldsymbol{\beta}})^2$, where $\widehat{\boldsymbol{\beta}}$ is the estimated coefficients vector by a method and $\mathbf{x}$ is an independent test point. The second measure, #S, is the number of selected variables in the final model by a method in a simulation. The third one, FN, measures the number of missed true variables by a method in a simulation.

In the calculation of PE, an independent test sample of size 10,000 was generated. All four methods had median FN = 0. Table 2 and Figure 3 summarize the comparison results given by PE and #S.

TABLE 1
*Success percentages of $L_1$, $\rho_a$ with different a and optimal SICA (with a tuned) in recovering $\boldsymbol{\beta}_0$ under three levels of correlation r in Simulation* 1*, where* $(s, n, p) = (7, 35, 1000)$

| | | SIRS | | | LLA | |
|---|---|---|---|---|---|---|
| **Methods** | $r = 0$ | $r = 0.2$ | $r = 0.5$ | $r = 0$ | $r = 0.2$ | $r = 0.5$ |
| $L_1$ | 11% | 8% | 4% | | | |
| $\rho_5$ | 19% | 12% | 8% | 28% | 26% | 20% |
| $\rho_2$ | 50% | 37% | 32% | 53% | 40% | 34% |
| $\rho_1$ | 63% | 58% | 51% | 57% | 46% | 43% |
| $\rho_{0.6}$ | 70% | 64% | 57% | 54% | 43% | 42% |
| $\rho_{0.4}$ | 74% | 69% | 62% | 51% | 38% | 38% |
| $\rho_{0.3}$ | 75% | 68% | 60% | 37% | 36% | 31% |
| $\rho_{0.2}$ | 76% | 70% | 61% | 27% | 24% | 17% |
| $\rho_{0.1}$ | 77% | 68% | 63% | 7% | 10% | 8% |
| $\rho_{0.05}$ | 71% | 68% | 61% | 2% | 3% | 2% |
| $\rho_0$ | 53% | 51% | 48% | — | — | — |
| Optimal SICA | 80% | 76% | 71% | 58% | 49% | 47% |

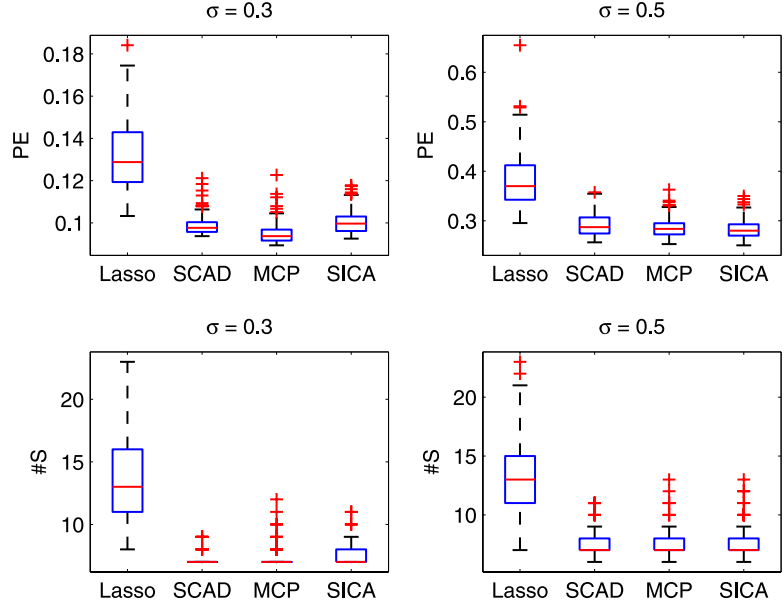

Fig. 3. *Boxplots of PE and #S over* 100 *simulations for all methods in Simulation* 2, *where* $p = 50$ *and the rows of* **X** *are i.i.d. copies from* $N(\mathbf{0}, \Sigma_0)$. *The x-axis represents different methods. Top panel is for PE and bottom panel is for #S.*

6.3. *Simulation* 3. The setting of this example is the same as that of Simulation 2, except that $(n, p) = (100, 600)$ and $\sigma = 0.1, 0.3$. Since $p$ is larger than $n$, BIC breaks down in the tuning of $\lambda$ and $a$. Thus we used five-fold cross-validation based on prediction error to select the tuning parameters. All four methods had median FN = 0. Table 3 and Figure 4 summarize the comparison results given by PE and #S. The boxplots of lasso are truncated to make it easier to view.

6.4. *Real data analysis.* In this example, we apply SICA to the diabetes dataset, which was studied by Efron et al. (2004). This dataset contains

TABLE 2
*Medians of PE and #S over* 100 *simulations for all methods in Simulation* 2, *where* $p = 50$ *and the rows of* **X** *are i.i.d. copies from* $N(\mathbf{0}, \Sigma_0)$

|  | Measures | Lasso | SCAD | MCP | SICA |
|---|---|---|---|---|---|
| $\sigma = 0.3$ | PE ($\times 10^{-2}$) | 12.88 | 9.77 | 9.38 | 9.97 |
|  | #S | 13 | 7 | 7 | 7 |
| $\sigma = 0.5$ | PE ($\times 10^{-1}$) | 3.70 | 2.87 | 2.84 | 2.80 |
|  | #S | 13 | 7 | 7 | 7 |

10 baseline variables: age (`age`), sex (`sex`), body mass index (`bmi`), average blood pressure (`bp`) and 6 blood serum measurements (`tc`, `ldl`, `hdl`, `tch`, `ltg`, `glu`) for $n = 442$ diabetes patients, as well as the response variable, a quantitative measure of disease progression one year after baseline. We implemented lasso, SCAD, MCP and SICA. Five-fold cross-validation was used to select the tuning parameters. All four methods excluded variable `hdl` in their final models. Lasso selected the remaining 9 variables, whereas SCAD, MCP and SICA all selected the same 6 variables. The estimated coefficients by different methods are shown in Table 4. For a comparison, we also included the adjusted $R^2$ and average prediction error in five-fold cross-validation for each method in Table 4. We see that SCAD, MCP and SICA performed similarly on this real dataset, while lasso produced a different model.

## 7. Proofs.

7.1. *Proof of Proposition 1.* Let $\mathbf{X} = \mathbf{UDV}$ be a singular value decomposition of $\mathbf{X}$, where $\mathbf{U}$ and $\mathbf{V}$ are, respectively, $n \times n$ and $p \times p$ orthogonal matrices, $\mathbf{D}$ is an $n \times p$ matrix with its first $k$ diagonal elements being $d_1, \ldots, d_k \neq 0$ and all other elements being zero, and $k = \mathrm{rank}(\mathbf{X}) \leq n \wedge p$. Then we have

$$
\begin{aligned}
\boldsymbol{\beta} \in \mathcal{A} &\iff \mathbf{y} = \mathbf{X}\boldsymbol{\beta} \iff \mathbf{U}^T \mathbf{y} = \mathbf{U}^T \mathbf{X} \boldsymbol{\beta} = \mathbf{D}\mathbf{V}\boldsymbol{\beta} = \mathbf{D}\widetilde{\boldsymbol{\beta}} \\
&\iff \widetilde{\beta}_i = d_i^{-1} w_i, \qquad i = 1, \ldots, k,
\end{aligned}
\tag{44}
$$

where $\widetilde{\boldsymbol{\beta}} = (\widetilde{\beta}_1, \ldots, \widetilde{\beta}_p)^T = \mathbf{V}\boldsymbol{\beta}$ and $\mathbf{U}^T \mathbf{y} = (w_1, \ldots, w_n)^T$. Since $\mathbf{V}$ is orthogonal, $\|\widetilde{\boldsymbol{\beta}}\|_2 = \|\boldsymbol{\beta}\|_2$ always holds. Thus it follows from (44) that

$$\boldsymbol{\beta}_2 = \arg\min_{\boldsymbol{\beta} \in \mathcal{A}} \|\boldsymbol{\beta}\|_2 = \mathbf{V}^T \widetilde{\boldsymbol{\beta}}_0,$$

where $\widetilde{\boldsymbol{\beta}}_0 = (d_1^{-1} w_1, \ldots, d_k^{-1} w_k, 0, \ldots, 0)^T$. It remains to show $\mathbf{V}^T \widetilde{\boldsymbol{\beta}}_0 = (\mathbf{X}^T \times \mathbf{X})^+ \mathbf{X}^T \mathbf{y}$. By the above singular value decomposition of $\mathbf{X}$, we have $(\mathbf{X}^T \mathbf{X})^+ = \mathbf{V}^T \mathrm{diag}(\mathbf{a}) \mathbf{V}$, where $\mathbf{a} = (d_1^{-2}, \ldots, d_k^{-2}, 0, \ldots, 0)^T$. Therefore, it is immediate to see that

$$\mathrm{diag}(\mathbf{a}) \mathbf{D}^T \mathbf{U}^T \mathbf{y} = (d_1^{-1} w_1, \ldots, d_k^{-1} w_k, 0, \ldots, 0)^T = \widetilde{\boldsymbol{\beta}}_0$$

and thus

$$(\mathbf{X}^T \mathbf{X})^+ \mathbf{X}^T \mathbf{y} = \mathbf{V}^T \mathrm{diag}(\mathbf{a}) \mathbf{V} \mathbf{V}^T \mathbf{D}^T \mathbf{U}^T \mathbf{y} = \mathbf{V}^T \widetilde{\boldsymbol{\beta}}_0.$$

TABLE 3
*Medians of PE and #S over 100 simulations for all methods in Simulation 3, where $p = 600$ and the rows of $\mathbf{X}$ are i.i.d. copies from $N(\mathbf{0}, \Sigma_0)$*

|                | Measures            | Lasso | SCAD | MCP  | SICA |
|----------------|---------------------|-------|------|------|------|
| $\sigma = 0.1$ | PE ($\times 10^{-2}$) | 2.51  | 1.07 | 1.07 | 1.10 |
|                | #S                  | 71.5  | 7    | 7    | 7    |
| $\sigma = 0.3$ | PE ($\times 10^{-2}$) | 21.15 | 9.92 | 9.89 | 9.90 |
|                | #S                  | 69.5  | 16   | 10   | 8    |

TABLE 4
*Model coefficients obtained by all methods on the diabetes dataset, and their adjusted $R^2$ [$R^2$(adj)] and average prediction errors (APE) based on five-fold cross-validation*

| Methods | age  | sex    | bmi   | bp    | tc        | ldl   |
|---------|------|--------|-------|-------|-----------|-------|
| Lasso   | −6.4 | −235.9 | 521.8 | 321.0 | −568.6    | 301.6 |
| SCAD    | 0    | −226.2 | 529.9 | 327.1 | −757.6    | 538.3 |
| MCP     | 0    | −226.3 | 529.9 | 327.2 | −757.7    | 538.4 |
| SICA    | 0    | −219.5 | 531.7 | 323.3 | −743.1    | 525.0 |
|         | hdl  | tch    | ltg   | glu   | $R^2$(adj) | APE   |
| Lasso   | 0    | 143.9  | 669.6 | 66.8  | 50.73%    | 2956.9 |
| SCAD    | 0    | 0      | 804.1 | 0     | 50.82%    | 2939.5 |
| MCP     | 0    | 0      | 804.1 | 0     | 50.82%    | 2939.1 |
| SICA    | 0    | 0      | 800.2 | 0     | 50.82%    | 2935.8 |

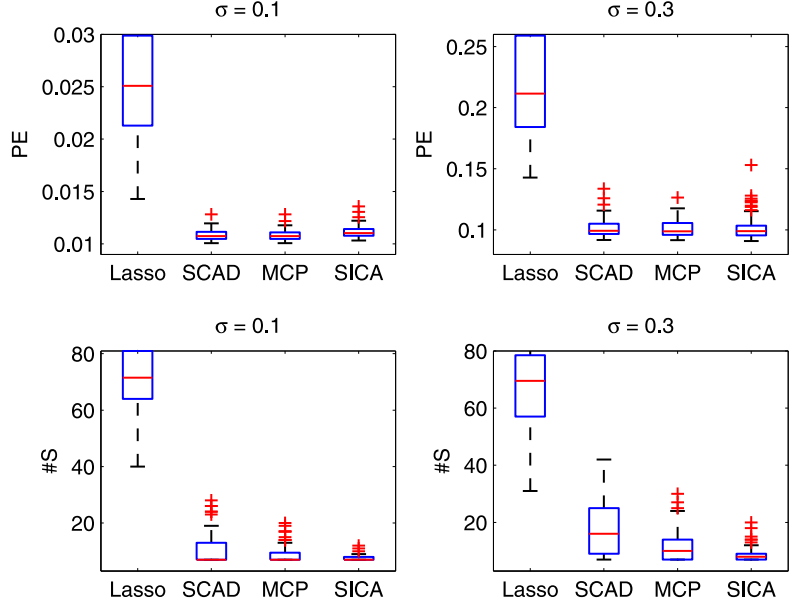

Fig. 4. *Boxplots of PE and #S over* 100 *simulations for all methods in Simulation* 3, *where* $p = 600$ *and the rows of* $\mathbf{X}$ *are i.i.d. copies from* $N(\mathbf{0}, \Sigma_0)$. *The x-axis represents different methods. Top panel is for PE and bottom panel is for #S.*

This completes the proof.

7.2. *Proof of Theorem 1.* Fix an arbitrary $\widehat{\boldsymbol{\beta}}^\lambda = (\widehat{\beta}_1^\lambda, \ldots, \widehat{\beta}_p^\lambda)^T \in \mathbf{R}^p$ and let $\widehat{\mathfrak{M}}_\lambda = \mathrm{supp}(\widehat{\boldsymbol{\beta}}^\lambda)$. It follows from the classical optimization theory by taking differentiation that if $\widehat{\boldsymbol{\beta}}^\lambda$ is a local minimizer of the regularized least squares problem (4) with $\lambda_n = \lambda$, there exists some $\mathbf{v} = (v_1, \ldots, v_p)^T \in \mathbf{R}^p$ such that

$$(45) \qquad \mathbf{X}^T \mathbf{X} \widehat{\boldsymbol{\beta}}^\lambda - \mathbf{X}^T \mathbf{y} + \Lambda_n \lambda \mathbf{v} = \mathbf{0},$$

where for $j \in \widehat{\mathfrak{M}}_\lambda$, $v_j = \bar{\rho}(\widehat{\beta}_j^\lambda)$ and for $j \in \widehat{\mathfrak{M}}_\lambda^c$, $v_j \in [-\rho'(0+), \rho'(0+)]$. Moreover, since $\widehat{\boldsymbol{\beta}}^\lambda$ is also a local minimizer of (4) constrained on the $\|\widehat{\boldsymbol{\beta}}^\lambda\|_0$-dimensional subspace $\{\boldsymbol{\beta} \in \mathbf{R}^p : \boldsymbol{\beta}_{\widehat{\mathfrak{M}}_\lambda^c} = \mathbf{0}\}$ of $\mathbf{R}^p$, it is easy to show that

$$(46) \qquad \lambda_{\min}(\mathbf{Q}) \geq \Lambda_n \lambda \kappa(\rho; \widehat{\boldsymbol{\beta}}_{\widehat{\mathfrak{M}}_\lambda}^\lambda),$$

where $\mathbf{Q} = \mathbf{X}_{\widehat{\mathfrak{M}}_\lambda}^T \mathbf{X}_{\widehat{\mathfrak{M}}_\lambda}$ and $\kappa(\rho; \widehat{\boldsymbol{\beta}}_{\widehat{\mathfrak{M}}_\lambda}^\lambda)$ is given by (8). We will see below that slightly strengthening the necessary condition (45) and (46) provides a sufficient condition on the strict local minimizer of (4).

Since $\mathbf{Q} = \mathbf{X}_{\widehat{\mathfrak{M}}_\lambda}^T \mathbf{X}_{\widehat{\mathfrak{M}}_\lambda}$ is nonsingular, $\widehat{\boldsymbol{\beta}}_{\widehat{\mathfrak{M}}_\lambda^c}^\lambda = \mathbf{0}$, $\mathbf{v}_{\widehat{\mathfrak{M}}_\lambda} = \bar{\rho}(\widehat{\boldsymbol{\beta}}_{\widehat{\mathfrak{M}}_\lambda}^\lambda)$ and $\|\mathbf{v}_{\widehat{\mathfrak{M}}_\lambda^c}\|_\infty \leq \rho'(0+)$, (45) can, equivalently, be rewritten as

$$(47) \qquad \widehat{\boldsymbol{\beta}}_{\widehat{\mathfrak{M}}_\lambda}^\lambda = \mathbf{Q}^{-1} \mathbf{X}_{\widehat{\mathfrak{M}}_\lambda}^T \mathbf{y} - \Lambda_n \lambda \mathbf{Q}^{-1} \bar{\rho}(\widehat{\boldsymbol{\beta}}_{\widehat{\mathfrak{M}}_\lambda}^\lambda),$$

$$(48) \qquad \|\mathbf{z}_{\widehat{\mathfrak{M}}_\lambda^c}\|_\infty \leq \rho'(0+),$$

where $\mathbf{z} = (\Lambda_n \lambda)^{-1} \mathbf{X}^T (\mathbf{y} - \mathbf{X} \widehat{\boldsymbol{\beta}}^\lambda)$. Now we strengthen inequality (48) to strict inequality (19) and make an additional assumption (20). We will show that (18), (19) and (20) imply that $\widehat{\boldsymbol{\beta}}^\lambda$ is a strict local minimizer of (4).

We first constrain the regularized least squares problem (4) on the $\|\widehat{\boldsymbol{\beta}}^\lambda\|_0$-dimensional subspace $\mathcal{B} = \{\boldsymbol{\beta} \in \mathbf{R}^p : \boldsymbol{\beta}_{\widehat{\mathfrak{M}}_\lambda^c} = \mathbf{0}\}$ of $\mathbf{R}^p$. It follows easily from condition (20), the continuity of $\rho'(t)$ in Condition 1, and the definition of $\kappa(\rho; \widehat{\boldsymbol{\beta}}_{\widehat{\mathfrak{M}}_\lambda}^\lambda)$ in (8) that the objective function in (4), $\ell(\boldsymbol{\beta}) \equiv 2^{-1} \|\mathbf{y} - \mathbf{X} \boldsymbol{\beta}\|_2^2 + \Lambda_n \sum_{j=1}^p p_\lambda(|\beta_j|)$, is strictly convex in a ball $\mathcal{N}_0$ in the subspace $\mathcal{B}$ centered at $\widehat{\boldsymbol{\beta}}^\lambda$. This along with (18) immediately entails that $\widehat{\boldsymbol{\beta}}^\lambda$, as a critical point of $\ell(\cdot)$ in $\mathcal{B}$, is the unique minimizer of $\ell(\cdot)$ in the neighborhood $\mathcal{N}_0$. Thus we have shown that $\widehat{\boldsymbol{\beta}}^\lambda$ is a strict local minimizer of $\ell(\cdot)$ in the subspace $\mathcal{B}$.

It remains to prove that the sparse vector $\widehat{\boldsymbol{\beta}}^\lambda$ is indeed a strict local minimizer of $\ell(\cdot)$ on the whole space $\mathbf{R}^p$. To show this, we will use condition (19). Take a sufficiently small ball $\mathcal{N}_1$ in $\mathbf{R}^p$ centered at $\widehat{\boldsymbol{\beta}}^\lambda$ such that $\mathcal{N}_1 \cap \mathcal{B} \subset \mathcal{N}_0$. Fix an arbitrary $\boldsymbol{\gamma}_1 \in \mathcal{N}_1 \setminus \mathcal{N}_0$, we will show that $\ell(\boldsymbol{\gamma}_1) > \ell(\widehat{\boldsymbol{\beta}}^\lambda)$.

Let $\boldsymbol{\gamma}_2$ be the projection of $\boldsymbol{\gamma}_1$ onto the subspace $\mathcal{B}$. Then it follows from $\mathcal{N}_1 \cap \mathcal{B} \subset \mathcal{N}_0$ and the definitions of $\mathcal{B}$, $\mathcal{N}_0$ and $\mathcal{N}_1$ that $\boldsymbol{\gamma}_2 \in \mathcal{N}_0 \cap \mathcal{N}_1$, which entails that $\ell(\boldsymbol{\gamma}_2) > \ell(\widehat{\boldsymbol{\beta}}^\lambda)$ if $\boldsymbol{\gamma}_2 \neq \widehat{\boldsymbol{\beta}}^\lambda$ by the strict convexity of $\ell(\cdot)$ in the neighborhood $\mathcal{N}_0$. We see that to prove $\ell(\boldsymbol{\gamma}_1) > \ell(\widehat{\boldsymbol{\beta}}^\lambda)$, it suffices to show that $\ell(\boldsymbol{\gamma}_1) > \ell(\boldsymbol{\gamma}_2)$.

It follows from the concavity of $\rho$ in Condition 1 that $\rho'(t)$ is decreasing in $t \in [0, \infty)$. By condition (19) and the continuity of $\rho'(t)$ in Condition 1, appropriately shrinking the radius of the ball $\mathcal{N}_1$ gives that there exists some $\delta \in (0, \infty)$ such that $\rho'(\delta) \in (0, \rho'(0+)]$ and for any $\boldsymbol{\beta} \in \mathcal{N}_1$,

$$\|\mathbf{w}_{\widehat{\mathfrak{M}}_\lambda^c}\|_\infty < \rho'(\delta), \tag{49}$$

where $\mathbf{w} = (\Lambda_n \lambda)^{-1} \mathbf{X}^T (\mathbf{y} - \mathbf{X}\boldsymbol{\beta})$. We further shrink the radius of the ball $\mathcal{N}_1$ to less than $\delta$. By the mean-value theorem, we have

$$\ell(\boldsymbol{\gamma}_1) = \ell(\boldsymbol{\gamma}_2) + \nabla^T \ell(\boldsymbol{\gamma}_0)(\boldsymbol{\gamma}_1 - \boldsymbol{\gamma}_2), \tag{50}$$

where $\boldsymbol{\gamma}_0$ lies on the line segment joining $\boldsymbol{\gamma}_2$ and $\boldsymbol{\gamma}_1$ and $\boldsymbol{\gamma}_0 \neq \boldsymbol{\gamma}_2$. Since $\boldsymbol{\gamma}_1, \boldsymbol{\gamma}_2 \in \mathcal{N}_1$ and $\mathcal{N}_1$ is a ball centered at $\widehat{\boldsymbol{\beta}}^\lambda$ with radius less than $\delta$, we have $\boldsymbol{\gamma}_0 = (\gamma_{0,1}, \ldots, \gamma_{0,p})^T \in \mathcal{N}_1$ and $|\gamma_{0,j}| < \delta$ for any $j \in \widehat{\mathfrak{M}}_\lambda^c$. Note that by $\boldsymbol{\gamma}_1 \in \mathcal{N}_1 \setminus \mathcal{N}_0$, we have

$$(\boldsymbol{\gamma}_1 - \boldsymbol{\gamma}_2)_{\widehat{\mathfrak{M}}_\lambda} = \mathbf{0} \quad \text{and} \quad (\boldsymbol{\gamma}_1 - \boldsymbol{\gamma}_2)_{\widehat{\mathfrak{M}}_\lambda^c} \neq \mathbf{0}.$$

Let $S = \text{supp}(\boldsymbol{\gamma}_1) \setminus \widehat{\mathfrak{M}}_\lambda \neq \varnothing$ and $\boldsymbol{\gamma}_1 = (\gamma_{1,1}, \ldots, \gamma_{1,p})^T$. It is easy to see that $\text{sgn}(\gamma_{0,j}) = \text{sgn}(\gamma_{1,j})$ for any $j \in \widehat{\mathfrak{M}}_\lambda^c$. Since $\boldsymbol{\gamma}_0 \in \mathcal{N}_1$, by (49), (50) and $\|\boldsymbol{\gamma}_{0,\widehat{\mathfrak{M}}_\lambda^c}\|_\infty < \delta$, we have

$$\begin{aligned}
\ell(\boldsymbol{\gamma}_1) - \ell(\boldsymbol{\gamma}_2) &= \sum_{j \in S} \frac{\partial \ell(\boldsymbol{\gamma}_0)}{\partial \beta_j} \gamma_j = [\mathbf{X}_S^T \mathbf{X} \boldsymbol{\gamma}_0 - \mathbf{X}_S^T \mathbf{y} + \Lambda_n \lambda \bar{\rho}(\boldsymbol{\gamma}_{0,S})]^T \boldsymbol{\gamma}_{1,S} \\
&= -\Lambda_n \lambda [(\Lambda_n \lambda)^{-1} \mathbf{X}_S^T (\mathbf{y} - \mathbf{X} \boldsymbol{\gamma}_0)]^T \boldsymbol{\gamma}_{1,S} \\
&\quad + \Lambda_n \lambda \sum_{j \in S} \text{sgn}(\gamma_{0,j}) \rho'(|\gamma_{0,j}|) \gamma_{1,j} \\
&> -\Lambda_n \lambda \rho'(\delta) \|\boldsymbol{\gamma}_{1,S}\|_1 + \Lambda_n \lambda \sum_{j \in S} \rho'(|\gamma_{0,j}|) |\gamma_{1,j}| \\
&\geq -\Lambda_n \lambda \rho'(\delta) \|\boldsymbol{\gamma}_{1,S}\|_1 + \Lambda_n \lambda \sum_{j \in S} \rho'(\delta) |\gamma_{1,j}| = 0,
\end{aligned}$$

where we used the fact that $\rho'(|\gamma_{0,j}|) \geq \rho'(\delta)$ since $\rho'(t)$ is decreasing in $t \in [0, \infty)$ and $|\gamma_{0,j}| < \delta$ for any $j \in S$. This shows that $\ell(\boldsymbol{\gamma}_1) > \ell(\boldsymbol{\gamma}_2)$, which concludes the proof.

7.3. *Proof of Theorem 2.* By (2), we have $\mathbf{y} = \mathbf{X}\boldsymbol{\beta}_0$. As mentioned before, to study the $\rho$-regularization problem (3), we consider the related $\rho$-regularized least squares problem (4) with $\mathbf{y} = \mathbf{X}\boldsymbol{\beta}_0$, $\Lambda_n = 1$ and $p_\lambda(t) = \lambda\rho(t)$, $t \in [0, \infty)$. We will construct a sequence of strict local minimizers $\widehat{\boldsymbol{\beta}}^\lambda = (\widehat{\beta}_1^\lambda, \ldots, \widehat{\beta}_p^\lambda)^T \in \mathbf{R}^p$ of (4) for a sequence of $\lambda \in (0, \infty)$ by using Theorem 1 and show that $\lim_{\lambda \to 0+} \widehat{\boldsymbol{\beta}}^\lambda = \boldsymbol{\beta}_0$. Moreover, with some careful analysis we show that the limit $\boldsymbol{\beta}_0$ is indeed a local minimizer of (3).

Fix an arbitrary $\lambda \in (0, \infty)$ and $\widehat{\boldsymbol{\beta}}^\lambda = (\widehat{\beta}_1^\lambda, \ldots, \widehat{\beta}_p^\lambda)^T \in \mathbf{R}^p$. By Theorem 1, $\widehat{\boldsymbol{\beta}}^\lambda$ will be a strict local minimizer of (4) as long as it satisfies conditions (18)–(20). We prove the existence of such a solution $\widehat{\boldsymbol{\beta}}^\lambda$ when $\widehat{\mathfrak{M}}_\lambda = \mathrm{supp}(\widehat{\boldsymbol{\beta}}^\lambda) = \mathfrak{M}_0$. Since $\mathbf{y} = \mathbf{X}\boldsymbol{\beta}_0$ and $\kappa(\rho; \widehat{\boldsymbol{\beta}}^\lambda_{\widehat{\mathfrak{M}}_\lambda}) \leq \kappa(\rho)$ in view of (7) and (8), we can rewrite (18) and (19) and strengthen (20) as

$$\widehat{\boldsymbol{\beta}}^\lambda_{\mathfrak{M}_0} = \boldsymbol{\beta}_{0,\mathfrak{M}_0} - \lambda \mathbf{Q}^{-1} \bar{\rho}(\widehat{\boldsymbol{\beta}}^\lambda_{\mathfrak{M}_0}), \tag{51}$$

$$\|\mathbf{z}_{\mathfrak{M}_0^c}\|_\infty < \rho'(0+), \tag{52}$$

$$\lambda_{\min}(\mathbf{Q}) > \lambda \kappa(\rho), \tag{53}$$

where $\mathbf{z} = \mathbf{X}^T \mathbf{X}_{\mathfrak{M}_0} \mathbf{Q}^{-1} \bar{\rho}(\widehat{\boldsymbol{\beta}}^\lambda_{\mathfrak{M}_0})$ and $\mathbf{Q} = \mathbf{X}_{\mathfrak{M}_0}^T \mathbf{X}_{\mathfrak{M}_0}$. Since $\mathbf{Q}$ is nonsingular and $\kappa(\rho) < \infty$ by assumption, all $\lambda \in (0, \lambda_0)$ with $\lambda_0 = \lambda_{\min}(\mathbf{Q})/\kappa(\rho)$ automatically satisfy condition (53). We now consider $\lambda \in (0, \lambda_0)$. Assume that there exists some $\epsilon \in (0, \min_{j \in \mathfrak{M}_0} |\beta_{0,j}|)$ such that (22) holds. We will show that there exists a solution $\widehat{\boldsymbol{\beta}}^\lambda_{\mathfrak{M}_0} \in \mathcal{V}_\epsilon$ to (51) and (52), where $\mathcal{V}_\epsilon = \prod_{j \in \mathfrak{M}_0} \{t : |t - \beta_{0,j}| \leq \epsilon\}$. Note that condition (22) guarantees (52). So it remains to prove the existence of $\widehat{\boldsymbol{\beta}}^\lambda_{\mathfrak{M}_0} \in \mathcal{V}_\epsilon$ to (51).

Note that $\rho'(t)$ is decreasing in $t \in [0, \infty)$ and thus $\|\bar{\rho}(\widehat{\boldsymbol{\beta}}^\lambda_{\mathfrak{M}_0})\|_\infty \leq \rho'(0+)$. Let $s = \|\boldsymbol{\beta}_0\|_0$,

$$h = \max_{\mathbf{v} \in \mathbf{R}^s, \|\mathbf{v}\|_\infty \leq \rho'(0+)} \|\mathbf{Q}^{-1}\mathbf{v}\|_\infty$$

and $\lambda_1 = \lambda_0 \vee (\epsilon/h)$. We now consider $\lambda \in (0, \lambda_1)$. Then for any $\boldsymbol{\gamma} \in \mathcal{V}_\epsilon$, we have

$$\|\lambda \mathbf{Q}^{-1} \bar{\rho}(\boldsymbol{\gamma})\|_\infty \leq \lambda_1 \|\mathbf{Q}^{-1} \bar{\rho}(\boldsymbol{\gamma})\|_\infty \leq \lambda_1 h \leq \epsilon.$$

Thus by the continuity of the vector-valued function $\boldsymbol{\Psi}(\boldsymbol{\gamma}) \equiv \boldsymbol{\gamma} - \boldsymbol{\beta}_{0,\mathfrak{M}_0} + \lambda \mathbf{Q}^{-1} \bar{\rho}(\boldsymbol{\gamma})$, an application of Miranda's existence theorem [see, e.g., Vrahatis (1989)] shows that (51) indeed has a solution $\widehat{\boldsymbol{\beta}}^\lambda_{\mathfrak{M}_0}$ in $\mathcal{V}_\epsilon$.

For any $\lambda \in (0, \lambda_1)$, we have shown that (4) has a strict local minimizer $\widehat{\boldsymbol{\beta}}^\lambda$ such that $\mathrm{supp}(\widehat{\boldsymbol{\beta}}^\lambda) = \mathfrak{M}_0$, $\widehat{\boldsymbol{\beta}}^\lambda_{\mathfrak{M}_0} \in \mathcal{V}_\epsilon$, and (51) holds. Note that by (51),

$$\|\widehat{\boldsymbol{\beta}}^\lambda - \boldsymbol{\beta}_0\|_\infty = \|\widehat{\boldsymbol{\beta}}^\lambda_{\mathfrak{M}_0} - \boldsymbol{\beta}_{0,\mathfrak{M}_0}\|_\infty = \lambda \|\mathbf{Q}^{-1} \bar{\rho}(\widehat{\boldsymbol{\beta}}^\lambda_{\mathfrak{M}_0})\|_\infty \leq \lambda h,$$

which entails that $\lim_{\lambda \to 0+} \widehat{\boldsymbol{\beta}}^\lambda$ exists and equals $\boldsymbol{\beta}_0$. It remains to prove that the limit $\boldsymbol{\beta}_0$ is indeed a local minimizer of (3).

In view of the choice of $\lambda$ and condition (22), it follows easily from (52), (53) and the proof of Theorem 1 that for each $\lambda \in (0, \lambda_1)$, $\widehat{\boldsymbol{\beta}}^\lambda$ is the strict minimizer of $\ell(\boldsymbol{\beta}) \equiv 2^{-1}\|\mathbf{y} - \mathbf{X}\boldsymbol{\beta}\|_2^2 + \lambda \rho(\boldsymbol{\beta})$ on some common neighborhood $\mathcal{C} = \{\boldsymbol{\beta} \in \mathbf{R}^p : \boldsymbol{\beta}_{\mathfrak{M}_0} \in \mathcal{V}_\epsilon, \|\boldsymbol{\beta}_{\mathfrak{M}_0^c}\|_\infty \le \delta\}$ of $\boldsymbol{\beta}_{n,0}$, for some $\delta \in (0, \infty)$ independent of $\lambda$. We will show that $\boldsymbol{\beta}_0$ is the minimizer of (3) on the neighborhood $\mathcal{N} = \mathcal{C} \cap \mathcal{A}$ of $\boldsymbol{\beta}_0$ in the subspace $\mathcal{A} = \{\boldsymbol{\beta} \in \mathbf{R}^p : \mathbf{y} = \mathbf{X}\boldsymbol{\beta}\}$. Fix an arbitrary $\boldsymbol{\gamma} \in \mathcal{N}$. Since $\boldsymbol{\gamma} \in \mathcal{A}$ and $\boldsymbol{\gamma}, \widehat{\boldsymbol{\beta}}^\lambda \in \mathcal{C}$ for each $\lambda \in (0, \lambda_1)$, it follows from (51), $\mathbf{y} = \mathbf{X}\boldsymbol{\beta}_0$ and $\mathbf{Q} = \mathbf{X}_{\mathfrak{M}_0}^T \mathbf{X}_{\mathfrak{M}_0}$ that for each $\lambda \in (0, \lambda_1)$,

$$\begin{aligned}(54)\quad \rho(\boldsymbol{\gamma}) &= \lambda^{-1}\ell(\boldsymbol{\gamma}) \ge \lambda^{-1}\ell(\widehat{\boldsymbol{\beta}}^\lambda) = \rho(\widehat{\boldsymbol{\beta}}^\lambda) + (2\lambda)^{-1}\|\mathbf{y} - \mathbf{X}\widehat{\boldsymbol{\beta}}^\lambda\|_2^2 \\ &= \rho(\widehat{\boldsymbol{\beta}}^\lambda) + 2^{-1}\lambda \bar{\rho}(\widehat{\boldsymbol{\beta}}_{\mathfrak{M}_0}^\lambda)^T \mathbf{Q}^{-1} \bar{\rho}(\widehat{\boldsymbol{\beta}}_{\mathfrak{M}_0}^\lambda).\end{aligned}$$

Note that by $\|\bar{\rho}(\widehat{\boldsymbol{\beta}}_{\mathfrak{M}_0}^\lambda)\|_\infty \le \rho'(0+)$, we have

$$0 \le \bar{\rho}(\widehat{\boldsymbol{\beta}}_{\mathfrak{M}_0}^\lambda)^T \mathbf{Q}^{-1} \bar{\rho}(\widehat{\boldsymbol{\beta}}_{\mathfrak{M}_0}^\lambda) \le \lambda_{\min}(\mathbf{Q})^{-1} \|\bar{\rho}(\widehat{\boldsymbol{\beta}}_{\mathfrak{M}_0}^\lambda)\|_2^2$$
$$\le \lambda_{\min}(\mathbf{Q})^{-1} s \rho'(0+)^2.$$

Thus by $\lim_{\lambda \to 0+} \widehat{\boldsymbol{\beta}}^\lambda = \boldsymbol{\beta}_0$ and the continuity of $\rho$, letting $\lambda \to 0+$ in (54) yields $\rho(\boldsymbol{\gamma}) \ge \rho(\boldsymbol{\beta}_0)$. This completes the proof.

7.4. *Proof of Theorem 3.* By (25) and (12), the optimal penalty $\rho_{a_{\mathrm{opt}}(\epsilon)}$ satisfies

$$\begin{aligned}(55)\quad \kappa(\rho_{a_{\mathrm{opt}}(\epsilon)}) &= \inf_{\rho_a \in \mathcal{P}_\epsilon} \kappa(\rho_a) = \inf_{\rho_a \in \mathcal{P}_\epsilon} 2(a^{-1} + a^{-2}) \\ &= \sup\{a \in (0, \infty] : \rho_a \in \mathcal{P}_\epsilon\}.\end{aligned}$$

Thus it follows from the definition of $\mathcal{P}_\epsilon$ in (24) that the optimal parameter $a_{\mathrm{opt}}(\epsilon)$ is the largest $a \in (0, \infty]$ such that the following condition holds:

$$\max_{j \in \mathfrak{M}_0^c} \max_{\mathbf{u} \in \mathcal{U}_\epsilon} |\langle \mathbf{x}_j, \mathbf{u}\rangle| \le \rho_a'(0+) = 1 + a^{-1},$$

where $\mathcal{U}_\epsilon = \{\mathbf{X}_{\mathfrak{M}_0} \mathbf{Q}^{-1} \bar{\rho}_a(\mathbf{v}) : \mathbf{v} \in \mathcal{V}_\epsilon\}$ and $\mathcal{V}_\epsilon = \prod_{j \in \mathfrak{M}_0} \{t : |t - \beta_{0,j}| \le \epsilon\}$. It is easy to see that for any $\epsilon \in (0, \min_{j \in \mathfrak{M}_0} |\beta_{0,j}|)$, we have $a_{\mathrm{opt}}(\epsilon) > 0$ since $1 + a^{-1} \to \infty$ and for each $t \ne 0$, $\bar{\rho}_a(t) = t^{-2} O(a)$ as $a \to 0+$. This proves part (a).

Part (b) follows from two simple facts. First, it is clear that $a_{\mathrm{opt}}(\epsilon) = \infty$ if and only if

$$(56)\quad \max_{j \in \mathfrak{M}_0^c} \max_{\mathbf{u} \in \mathcal{U}_\epsilon} |\langle \mathbf{x}_j, \mathbf{u}\rangle| \le 1,$$

where $\mathcal{U}_\epsilon = \{\mathbf{X}_{\mathfrak{M}_0}\mathbf{Q}^{-1}\mathrm{sgn}(\mathbf{v}) : \mathbf{v} \in \mathcal{V}_\epsilon\}$ and $\mathcal{V}_\epsilon = \prod_{j \in \mathfrak{M}_0}\{t : |t - \beta_{0,j}| < \epsilon\}$. Second, for any $\epsilon \in (0, \min_{j \in \mathfrak{M}_0}|\beta_{0,j}|)$, $\mathcal{U}_\epsilon$ contains a single point $\mathbf{u}_0 = \mathbf{X}_{\mathfrak{M}_0}\mathbf{Q}^{-1}\mathrm{sgn}(\boldsymbol{\beta}_{0,\mathfrak{M}_0})$.

7.5. *Proof of Theorem 4.* We will prove that under the given regularity conditions, there exists a solution $\widehat{\boldsymbol{\beta}}^{\lambda_n} \in \mathbf{R}^p$ to (31)–(33) with $\widehat{\mathfrak{M}}_{\lambda_n} = \mathrm{supp}(\widehat{\boldsymbol{\beta}}^{\lambda_n}) = \mathfrak{M}_0$. Consider events

$$\mathcal{E}_1 = \{\|\mathbf{X}_{\mathfrak{M}_0}^T\boldsymbol{\varepsilon}\|_\infty \leq u_n D_{1n}\sigma\} \quad \text{and} \quad \mathcal{E}_2 = \{\|\mathbf{X}_{\mathfrak{M}_0^c}^T\boldsymbol{\varepsilon}\|_\infty \leq u_n D_{2n}\sigma\},$$

where $D_{1n}$ and $D_{2n}$ are defined in Condition 4. Let $\boldsymbol{\xi} = (\xi_1, \ldots, \xi_p)^T = \mathbf{X}^T\boldsymbol{\varepsilon}$. It follows from the definitions of $D_{1n}$ and $D_{2n}$ that

$$\mathcal{F}_1 = \{|\xi_j| \leq u_n\|\mathbf{x}_j\|_2\sigma : j \in \mathfrak{M}_0\} \subset \mathcal{E}_1$$

and

$$\mathcal{F}_2 = \{|\xi_j| \leq u_n\|\mathbf{x}_j\|_2\sigma : j \in \mathfrak{M}_0^c\} \subset \mathcal{E}_2,$$

where $\mathbf{X} = (\mathbf{x}_1, \ldots, \mathbf{x}_p)$. Since $\boldsymbol{\varepsilon} \sim N(\mathbf{0}, \sigma^2 I_n)$ by Condition 3, we see that for each $j = 1, \ldots, p$, $\xi_j$ has a $N(0, \|\mathbf{x}_j\|_2^2\sigma^2)$ distribution. Thus an application of the classical standard Gaussian tail probability bound and Bonferroni's inequality gives

$$\begin{aligned}(57) \quad P(\mathcal{E}_1 \cap \mathcal{E}_2) &\geq P(\mathcal{F}_1 \cap \mathcal{F}_2) \geq 1 - [P(\mathcal{F}_1^c) + P(\mathcal{F}_2^c)] \\ &\geq 1 - [2sP(V > u_n) + 2(p-s)P(V > u_n)] \\ &\geq 1 - \frac{2}{\sqrt{\pi}}pu_n^{-1}e^{-u_n^2/2},\end{aligned}$$

where $s = \|\boldsymbol{\beta}_0\|_0$ and $V$ is a standard Gaussian random variable. Hereafter we condition on the event $\mathcal{E}_1 \cap \mathcal{E}_2$. Under this event, we will show the existence of a solution $\widehat{\boldsymbol{\beta}}^{\lambda_n} \in \mathbf{R}^p$ to (31)–(33) with $\mathrm{sgn}(\widehat{\boldsymbol{\beta}}^{\lambda_n}) = \mathrm{sgn}(\boldsymbol{\beta}_0)$ and $\|\widehat{\boldsymbol{\beta}}^{\lambda_n} - \boldsymbol{\beta}_0\|_\infty \leq (1-c_0)b_0$.

By Condition 4, we have $\underline{\lambda}_n \leq \overline{\lambda}_n$, where

$$(58) \quad \underline{\lambda}_n = \Lambda_n^{-1}\frac{(C_{2n}D_{1n} + D_{2n})u_n\sigma}{\rho'(0+) - C_{2n}\rho'(c_0 b_0)} \quad \text{and} \quad \overline{\lambda}_n = \frac{C_{1n}^{-1}(1-c_0)b_0 - u_n D_{1n}\sigma}{\Lambda_n \rho'(c_0 b_0; \overline{\lambda}_n)}.$$

Let $\lambda_n$ be in the interval $[\underline{\lambda}_n, \overline{\lambda}_n]$. Since $\mathbf{y} = \mathbf{X}\boldsymbol{\beta}_0 + \boldsymbol{\varepsilon}$ by (1), (31) and (32) with $\widehat{\mathfrak{M}}_{\lambda_n} = \mathfrak{M}_0$ becomes

$$(59) \qquad \widehat{\boldsymbol{\beta}}_{\mathfrak{M}_0}^{\lambda_n} = \boldsymbol{\beta}_{0,\mathfrak{M}_0} - \mathbf{v},$$

$$(60) \qquad \|\mathbf{z}\|_\infty \leq \rho'(0+),$$

where $\mathbf{v} = (\mathbf{X}_{\mathfrak{M}_0}^T \mathbf{X}_{\mathfrak{M}_0})^{-1}[\Lambda_n \lambda_n \bar{\rho}(\widehat{\boldsymbol{\beta}}_{\mathfrak{M}_0}^{\lambda_n}) - \mathbf{X}_{\mathfrak{M}_0}^T \boldsymbol{\varepsilon}]$ and

$$\mathbf{z} = \mathbf{X}_{\mathfrak{M}_0^c}^T \mathbf{X}_{\mathfrak{M}_0} (\mathbf{X}_{\mathfrak{M}_0}^T \mathbf{X}_{\mathfrak{M}_0})^{-1} \bar{\rho}(\widehat{\boldsymbol{\beta}}_{\mathfrak{M}_0}^{\lambda_n})$$
$$- (\Lambda_n \lambda_n)^{-1}[\mathbf{X}_{\mathfrak{M}_0^c}^T \mathbf{X}_{\mathfrak{M}_0} (\mathbf{X}_{\mathfrak{M}_0}^T \mathbf{X}_{\mathfrak{M}_0})^{-1} \mathbf{X}_{\mathfrak{M}_0}^T \boldsymbol{\varepsilon} - \mathbf{X}_{\mathfrak{M}_0^c}^T \boldsymbol{\varepsilon}].$$

We first prove that (59) has a solution in $\mathcal{N} = \{\boldsymbol{\gamma} \in \mathbf{R}^s : \|\boldsymbol{\gamma} - \boldsymbol{\beta}_{0,\mathfrak{M}_0}\|_\infty \leq (1-c_0)b_0\}$. Fix an arbitrary $\boldsymbol{\gamma} \in \mathcal{N}$. By the concavity of $\rho$ in Condition 1, $\rho'(t)$ is decreasing in $t \in [0,\infty)$ and thus $\|\bar{\rho}(\boldsymbol{\gamma})\|_\infty \leq \rho'(c_0 b_0)$. This along with $\|\mathbf{X}_{\mathfrak{M}_0}^T \boldsymbol{\varepsilon}\|_\infty \leq u_n D_{1n}\sigma$, (34) in Condition 2 and $\lambda_n \leq \bar{\lambda}_n$ yields

$$\|(\mathbf{X}_{\mathfrak{M}_0}^T \mathbf{X}_{\mathfrak{M}_0})^{-1}[\Lambda_n \lambda_n \bar{\rho}(\boldsymbol{\gamma}) - \mathbf{X}_{\mathfrak{M}_0}^T \boldsymbol{\varepsilon}]\|_\infty \leq (1-c_0)b_0,$$

since $\rho'(t;\lambda)$ is increasing in $\lambda \in (0,\infty)$ by Condition 1. Thus by the continuity of the vector-valued function $\boldsymbol{\Psi}(\boldsymbol{\gamma}) \equiv \boldsymbol{\gamma} - \boldsymbol{\beta}_{0,\mathfrak{M}_0} + (\mathbf{X}_{\mathfrak{M}_0}^T \mathbf{X}_{\mathfrak{M}_0})^{-1}[\Lambda_n \lambda_n \bar{\rho}(\boldsymbol{\gamma}) - \mathbf{X}_{\mathfrak{M}_0}^T \boldsymbol{\varepsilon}]$, an application of Miranda's existence theorem [see, e.g., Vrahatis (1989)] shows that (59) indeed has a solution $\widehat{\boldsymbol{\beta}}_{\mathfrak{M}_0}^{\lambda_n}$ in $\mathcal{N}$. It remains to check the inequality (60) for $\widehat{\boldsymbol{\beta}}_{\mathfrak{M}_0}^{\lambda_n} \in \mathcal{N}$. In fact it can be easily shown to hold by (35) in Condition 2 and $\|\mathbf{X}_{\mathfrak{M}_0^c}^T \boldsymbol{\varepsilon}\|_\infty \leq u_n D_{2n}\sigma$, for $\lambda_n = \underline{\lambda}_n$.

So far we have shown the existence of a solution $\widehat{\boldsymbol{\beta}}^{\lambda_n}$ with $\lambda_n = \underline{\lambda}_n$ to (31) and (32) with $\mathrm{sgn}(\widehat{\boldsymbol{\beta}}^{\lambda_n}) = \mathrm{sgn}(\boldsymbol{\beta}_0)$ and $\|\widehat{\boldsymbol{\beta}}^{\lambda_n} - \boldsymbol{\beta}_0\|_\infty \leq (1-c_0)b_0$ under the event $\mathcal{E}_1 \cap \mathcal{E}_2$. By (58), (59) and $C_{2n} \leq C\frac{\rho'(0+)}{\rho'(c_0 b_0)}$, letting $\lambda_n = \underline{\lambda}_n$ gives

$$\|\widehat{\boldsymbol{\beta}}_{\mathfrak{M}_0}^{\lambda_n} - \boldsymbol{\beta}_{0,\mathfrak{M}_0}\|_\infty = \|(\mathbf{X}_{\mathfrak{M}_0}^T \mathbf{X}_{\mathfrak{M}_0})^{-1}[\Lambda_n \lambda_n \bar{\rho}(\widehat{\boldsymbol{\beta}}_{\mathfrak{M}_0}^{\lambda_n}) - \mathbf{X}_{\mathfrak{M}_0}^T \boldsymbol{\varepsilon}]\|_\infty$$
$$\leq u_n \left[D_{1n} + \frac{\rho'(c_0 b_0)}{\rho'(0+)} D_{2n}\right] C_{1n}(1-C)^{-1}\sigma.$$

Note that condition (33) with $\lambda_n = \underline{\lambda}_n$ is guaranteed by (37) in Condition 4 since $C_{2n} \leq C\frac{\rho'(0+)}{\rho'(c_0 b_0)}$. These along with $\mathrm{sgn}(\widehat{\boldsymbol{\beta}}^{\lambda_n}) = \mathrm{sgn}(\boldsymbol{\beta}_0)$, (36) in Condition 4 and (57) prove parts (a) and (b). Note that $1 - \frac{2}{\sqrt{\pi}} p u_n^{-1} e^{-u_n^2/2} \to 1$ since $p = o(u_n e^{u_n^2/2})$. Thus $\|\widehat{\boldsymbol{\beta}}^{\lambda_n} - \boldsymbol{\beta}_0\|_2 = O_P(\sqrt{s} n^{-\gamma} u_n)$ follows from parts (a) and (b) and

$$\|\widehat{\boldsymbol{\beta}}_{\mathfrak{M}_0}^{\lambda_n} - \boldsymbol{\beta}_{0,\mathfrak{M}_0}\|_2 \leq \sqrt{s}\|\widehat{\boldsymbol{\beta}}_{\mathfrak{M}_0}^{\lambda_n} - \boldsymbol{\beta}_{0,\mathfrak{M}_0}\|_\infty.$$

This concludes the proof.

7.6. *Proof of Proposition 2.* Part (a) follows directly from the definition of $\boldsymbol{\beta}^{(\ell)} = \mathbf{v}(\boldsymbol{\beta}^{(\ell-1)})$. It is not hard to check that $\boldsymbol{\beta}_0$ is the unique minimizer of the weighted $L_2$-regularization problem

$$\min_{\boldsymbol{\beta} \in \mathcal{A}} \boldsymbol{\beta}^T \boldsymbol{\Gamma}(\boldsymbol{\beta}_0) \boldsymbol{\beta},$$

which entails immediately part (b). It remains to prove part (c). Clearly, it suffices to show that for any $\boldsymbol{\beta} \in \mathcal{A}$, $\|\boldsymbol{\beta}\|_0 \leq (n+1)/2$ implies $\boldsymbol{\beta} = \boldsymbol{\beta}_0$. We prove this by a contradiction argument. Suppose there exists some $\boldsymbol{\beta} \in \mathcal{A}$ with $\|\boldsymbol{\beta}\|_0 \leq (n+1)/2$ and $\boldsymbol{\beta} \neq \boldsymbol{\beta}_0$. Let $\boldsymbol{\gamma} = \boldsymbol{\beta} - \boldsymbol{\beta}_0$. Then we have $\boldsymbol{\gamma} \neq \mathbf{0}$ and $\mathbf{X}\boldsymbol{\gamma} = \mathbf{X}\boldsymbol{\beta} - \mathbf{X}\boldsymbol{\beta}_0 = \mathbf{0}$. But

$$\|\boldsymbol{\gamma}\|_0 \leq \|\boldsymbol{\beta}\|_0 + \|\boldsymbol{\beta}_0\|_0 < (n+1)/2 + (n+1)/2 = n+1 = \text{spark}(\mathbf{X}),$$

which contradicts the definition of spark.

**8. Discussion.** We have studied the properties of regularization methods in model selection and sparse recovery under the unified framework of regularized least squares with concave penalties. We have provided regularity conditions under which the regularized least squares estimator enjoys a nonasymptotic weak oracle property for model selection, where the dimensionality can be of exponential order. Our results generalize those of Fan and Li (2001) and Fan and Peng (2004) in the setting of regularized least squares. For sparse recovery, we have generalized a sufficient condition identified for the $L_1$ penalty to concave penalties, which ensures the $\rho/L_0$ equivalence. In particular, a family of penalties that give a smooth homotopy between $L_0$ and $L_1$ penalties have been considered for both problems. Numerical studies further endorse our theoretical results and the advantage of our new methods for model selection and sparse recovery.

It would be interesting to extend the results to regularization methods for the generalized linear models (GLMs) and more general models and loss functions. These problems are beyond the scope of this paper and will be interesting topics for future research.

**Acknowledgments.** We sincerely thank Professors Jianqing Fan and Jun S. Liu for constructive comments that led to improvement over an earlier version of this paper. We are deeply indebted to Professor Jianqing Fan for introducing us the topic of high-dimensional variable selection and for his encouragement and helpful discussions over the years. We gratefully acknowledge the helpful comments of the Associate Editor and referees that substantially improved the presentation of the paper.

## REFERENCES

ANTONIADIS, A. and FAN, J. (2001). Regularization of wavelets approximations (with discussion). *J. Amer. Statist. Assoc.* **96** 939–967. MR1946364

BICKEL, P. J. and LI, B. (2006). Regularization in statistics (with discussion). *Test* **15** 271–344. MR2273731

BICKEL, P. J., RITOV, Y. and TSYBAKOV, A. (2008). Simultaneous analysis of Lasso and Dantzig selector. *Ann. Statist.* To appear.

BREIMAN, L. (1995). Better subset regression using the nonnegative garrote. *Technometrics* **37** 373–384. MR1365720

Candes, E. J. and Tao, T. (2005). Decoding by linear programming. *IEEE Trans. Inform. Theory* **51** 4203–4215. MR2243152

Candes, E. J. and Tao, T. (2006). Near-optimal signal recovery from random projections: Universal encoding strategies? *IEEE Trans. Inform. Theory* **52** 5406–5425. MR2300700

Candes, E. J. and Tao, T. (2007). The Dantzig selector: Statistical estimation when $p$ is much larger than $n$ (with discussion). *Ann. Statist.* **35** 2313–2404. MR2382644

Candès, E. J., Wakin, M. B. and Boyd, S. P. (2008). Enhancing sparsity by reweighted $\ell_1$ minimization. *J. Fourier Anal. Appl.* **14** 877–905.

Chen, S., Donoho, D. and Saunders, M. (1999). Atomic decomposition by basis pursuit. *SIAM J. Sci. Comput.* **20** 33–61. MR1639094

Donoho, D. L. (2004). Neighborly polytopes and sparse solution of underdetermined linear equations. Technical report, Dept. Statistics, Stanford Univ.

Donoho, D. L. and Elad, M. (2003). Optimally sparse representation in general (nonorthogonal) dictionaries via $\ell_1$ minimization. *Proc. Natl. Acad. Sci. USA* **100** 2197–2202. MR1963681

Donoho, D., Elad, M. and Temlyakov, V. (2006). Stable recovery of sparse overcomplete representations in the presence of noise. *IEEE Trans. Inform. Theory* **52** 6–18. MR2237332

Donoho, D. L. and Johnstone, I. M. (1994). Ideal spatial adaptation by wavelet shrinkage. *Biometrika* **81** 425–455. MR1311089

Efron, B., Hastie, T., Johnstone, I. and Tibshirani, R. (2004). Least angle regression (with discussion). *Ann. Statist.* **32** 407–451. MR2060166

Fan, J. (1997). Comment on "Wavelets in statistics: A review" by A. Antoniadis. *J. Italian Statist. Assoc.* **6** 131–138.

Fan, J. and Fan, Y. (2008). High-dimensional classification using features annealed independence rules. *Ann. Statist.* **36** 2605–2637.

Fan, J. and Li, R. (2001). Variable selection via nonconcave penalized likelihood and its oracle properties. *J. Amer. Statist. Assoc.* **96** 1348–1360. MR1946581

Fan, J. and Li, R. (2006). Statistical challenges with high dimensionality: Feature selection in knowledge discovery. In *Proceedings of the International Congress of Mathematicians* (M. Sanz-Sole, J. Soria, J. L. Varona and J. Verdera, eds.) **3** 595–622. European Math. Soc. Publishing House, Zürich. MR2275698

Fan, J. and Lv, J. (2008). Sure independence screening for ultrahigh dimensional feature space (with discussion). *J. Roy. Statist. Soc. Ser. B* **70** 849–911.

Fan, J. and Peng, H. (2004). Nonconcave penalized likelihood with diverging number of parameters. *Ann. Statist.* **32** 928–961. MR2065194

Fang, K.-T. and Zhang, Y.-T. (1990). *Generalized Multivariate Analysis*. Springer, Berlin. MR1079542

Frank, I. E. and Friedman, J. H. (1993). A statistical view of some chemometrics regression tools (with discussion). *Technometrics* **35** 109–148.

Fuchs, J.-J. (2004). Recovery of exact sparse representations in the presence of noise. In *Proceedings of IEEE International Conference on Acoustics, Speech, and Signal Processing* 533–536. Montreal, QC.

Greenshtein, E. and Ritov, Y. (2004). Persistence in high-dimensional linear predictor selection and the virtue of overparametrization. *Bernoulli* **10** 971–988. MR2108039

Hunter, D. and Li, R. (2005). Variable selection using MM algorithms. *Ann. Statist.* **33** 1617–1642. MR2166557

James, G., Radchenko, P. and Lv, J. (2009). DASSO: Connections between the Dantzig selector and Lasso. *J. Roy. Statist. Soc. Ser. B* **71** 127–142.

LI, R. and LIANG, H. (2008). Variable selection in semiparametric regression modeling. *Ann. Statist.* **36** 261–286. MR2387971

LIU, Y. and WU, Y. (2007). Variable selection via a combination of the $L_0$ and $L_1$ penalties. *J. Comput. Graph. Statist.* **16** 782–798.

MEINSHAUSEN, N., ROCHA, G. and YU, B. (2007). Discussion: A tale of three cousins: Lasso, L2Boosting and Dantzig. *Ann. Statist.* **35** 2373–2384. MR2382649

NIKOLOVA, M. (2000). Local strong homogeneity of a regularized estimator. *SIAM J. Appl. Math.* **61** 633–658. MR1780806

TIBSHIRANI, R. (1996). Regression shrinkage and selection via the Lasso. *J. Roy. Statist. Soc. Ser. B* **58** 267–288. MR1379242

TROPP, J. A. (2006). Just relax: Convex programming methods for identifying sparse signals in noise. *IEEE Trans. Inform. Theory* **5** 1030–1051. MR2238069

VRAHATIS, M. N. (1989). A short proof and a generalization of Miranda's existence theorem. *Proc. Amer. Math. Soc.* **107** 701–703. MR0993760

WAINWRIGHT, M. J. (2006). Sharp thresholds for high-dimensional and noisy recovery of sparsity. Technical report, Dept. Statistics, Univ. California, Berkeley.

WANG, H., LI, R. and TSAI, C.-L. (2007). Tuning parameter selectors for the smoothly clipped absolute deviation method. *Biometrika* **94** 553–568. MR2410008

ZHANG, C.-H. (2007). Penalized linear unbiased selection. Technical report, Dept. Statistics, Rutgers Univ.

ZHAO, P. and YU, B. (2006). On model selection consistency of Lasso. *J. Mach. Learn. Res.* **7** 2541–2563. MR2274449

ZOU, H. (2006). The adaptive Lasso and its oracle properties. *J. Amer. Statist. Assoc.* **101** 1418–1429. MR2279469

ZOU, H. and LI, R. (2008). One-step sparse estimates in nonconcave penalized likelihood models (with discussion). *Ann. Statist.* **36** 1509–1566. MR2435443

INFORMATION AND OPERATIONS
MANAGEMENT DEPARTMENT
MARSHALL SCHOOL OF BUSINESS
UNIVERSITY OF SOUTHERN CALIFORNIA
LOS ANGELES, CALIFORNIA 90089
USA
E-MAIL: jinchilv@marshall.usc.edu
    fanyingy@marshall.usc.edu